\documentclass[11pt]{article}
\usepackage{amsmath}
\usepackage{amsfonts}
\usepackage{graphicx}
\usepackage{latexsym}
\usepackage{color}

\newcommand{\R}{\mathbb{R}}
\newcommand{\E}{\mathbb{E}}
\newcommand{\PP}{\mathbb{P}}

\newcommand{\HH}{\mathbb{H}}

\newcommand{\De}{\Delta}

\newcommand{\stab}{\stackrel{d_{st}}{\longrightarrow}}
\newcommand{\schw}{\stackrel{d}{\longrightarrow}}
\newcommand{\toop}{\stackrel{\PP}{\longrightarrow}}

\newcommand{\ucp}{\stackrel{\mbox{\tiny u.c.p.}}{\Longrightarrow}}

\newcommand{\bee}{\begin{equation}}
\newcommand{\eee}{\end{equation}}
\newcommand{\bea}{\begin{eqnarray}}
\newcommand{\eea}{\end{eqnarray}}
\newcommand{\bean}{\begin{eqnarray*}}
\newcommand{\eean}{\end{eqnarray*}}

\newcommand{\qed}{$\hfill\Box$}

\newtheorem{prop}{Proposition}[section]

\newtheorem{defi}[prop]{Definition}

\newtheorem{theo}[prop]{Theorem}
\newtheorem{rem}[prop]{Remark}

\begin{document}

\title{Ambit fields: survey and new challenges}
\author{
Mark Podolskij \thanks{Department of Mathematics, University of Aarhus,
Ny Munkegade 118, 8000 Aarhus C,
Denmark, Email: mpodolskij@creates.au.dk.}}

\date{\today}

\maketitle

\begin{abstract}
In this paper we present a survey on recent developments in the study of ambit fields
and point out some open problems. Ambit fields is a class of spatio-temporal stochastic processes,
which by its general structure constitutes a flexible model for dynamical structures in time and/or in space.
We will review their basic probabilistic properties,  main stochastic integration concepts
and recent limit theory for high frequency statistics of ambit fields.

\ \

\bigskip


\end{abstract}

\section{Introduction}
\label{Intro}
\setcounter{equation}{0}
\renewcommand{\theequation}{\thesection.\arabic{equation}}
The recent years have witnessed a strongly increasing interest in ambit stochastics.
Ambit fields\footnote{From Latin \textit{ambitus}: a sphere of influence} 
is a class of spatio-temporal stochastic processes that has been
originally introduced by Barndorff-Nielsen and Schmiegel in a series of papers \cite{BS07, BS08a, BS09} in
the context of turbulence modelling, but which found manifold applications in 
mathematical finance and biology among other sciences; see e.g. \cite{BBC, BJJS07}.

\begin{figure}[t]
\centering
%\iffigures
\includegraphics[width=0.5\textwidth]{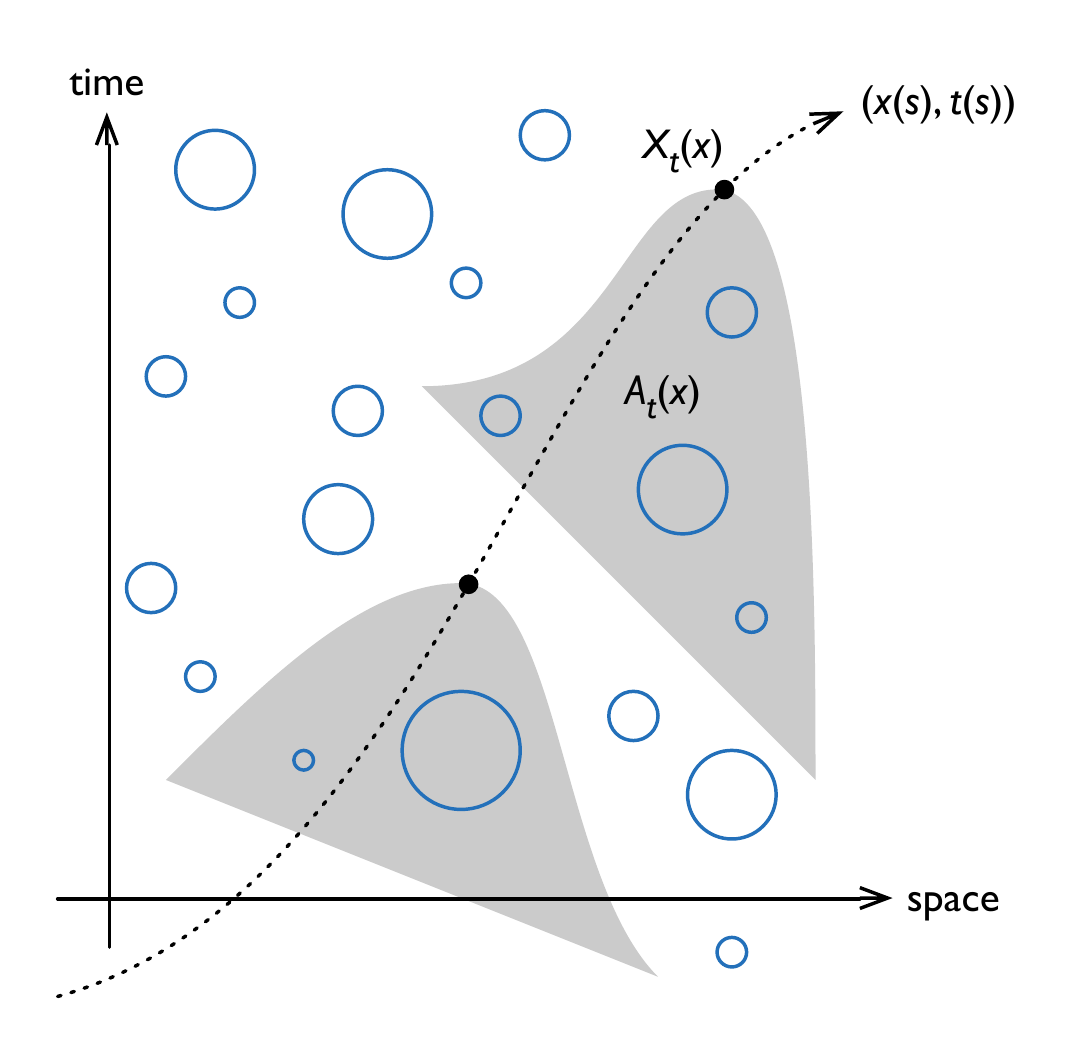} 
%\else
%\frame{\strut\ placeholder for ambit-figure.pdf } %\fi
\caption{A spatio-temporal ambit field. The value $X_t(x)$ of the
  field at the point $(x,t)$ marked by the nearby dot is defined
  through an integral over the corresponding ambit set $A_t(x)$ marked
  by the shaded region. The blue circles of varying sizes indicate the
  stochastic volatility/intermittency. By considering the field along
  the path $(x(s),t(s))$ in space-time an ambit process is
  obtained.}
\label{fig:af}
\end{figure}

Ambit processes describe the dynamics in a stochastically developing
field, for instance a turbulent wind field, along curves embedded in
such a field.  An illustration of the various concepts involved is
demonstrated by figure \ref{fig:af}. This shows a curve running with
time through space. To each point on the curve is associated a random
variable $X_{\theta }$\ which should be thought of as the value of an
observation at that point. This value is assumed to depend only on a
region of previous space-time, indicated by the solid line in the
figure, called the ambit set associated to the point in
question. Further, a key characteristic of the modelling framework,
which distinguishes ambit fields from other approaches is that
beyond the most basic kind of random input it also specifically
incorporates additional, often drastically changing, inputs referred
to as \textit{volatility} or \textit{intermittency}. This feature is illustrated by the
blue circles whose sizes indicate volatility effects of varying
degree.

In terms of mathematical formulae an ambit field is specified via
\begin{align} \label{ambfield}
  X_{t}(x)=\mu +\int_{A_{t}(x)}g(t,s,x,\xi )\sigma_{s}(\xi )L(\mathrm{d}s,\mathrm{d}
  \xi )+\int_{D_{t}(x)}q(t,s,x,\xi )a_{s}(\xi )\mathrm{d}s\mathrm{d}\xi
\end{align}
where $t$ denotes time while $x$\ gives the position in
space. Further, $A_{t}(x)$ and $D_{t}(x)$ are \textit{ambit sets}, $g$ and $q$
are deterministic weight functions, $\sigma$ represents the volatility or
intermittency field, $a$ is a drift field and $L$ denotes a \textit{L\'{e}vy
basis}. We recall that a L\'{e}vy basis $L=\{L(B):~ B\in \mathcal S\}$, where $\mathcal S$
is a $\delta$-ring of an arbitrary non-empty set $S$ such that there exists an increasing sequence
of sets $(S_n)\subset \mathcal S$ with $\cup_{n\in \mathbb N} S_n=S$, is an independently scattered random measure
with L\'evy-Khinchin representation
\begin{align} \label{charlevy}
\log \E[\exp(iu L(B))]= iu v_1(B) - \frac 12 u^2 v_2(B) + \int_{\R} \left( \exp(iuy) - 1- iuy 1_{[-1,1]}(y) \right) v_3(dy,B).
\end{align}
Here $v_1$ is a signed measure on $\mathcal S$,  $v_2$ is a measure on $\mathcal S$ and $v_3(\cdot,\cdot)$ is a generalised
L\'evy measure on $\R\times \mathcal S$ (see e.g. \cite{RR} for details).
In the turbulence framework the stochastic field $(X_{t}(x))_{t\geq 0, x\in \R^3}$ describes the velocity
of a turbulent flow at time $t$ and position $x$, while the ambit sets $A_{t}(x), D_{t}(x)$ are typically
bounded. 

Triggered by the wide applicability of ambit fields there is an
increasing need in understanding their mathematical properties. The
main mathematical topics of interest are: stochastic analysis for
ambit fields, modelling of volatility/intermittency, asymptotic theory
for ambit fields, numerical approximations of ambit fields, and asymptotic statistics.
The aim of this paper is to give an overview of recent studies and to point out the main challenges 
in future. 
We will see that most of the mathematical problems already arise in a null spatial 
setting, and therefore parts of the presented theory will deal with this case only. However,
possible extensions to spatial dimension will be discussed. At this stage we would like to mention
another review on ambit fields \cite{BBC2}, which focuses on integration theory, various examples of ambit
fields and their applications in physics and finance. Numerical schemes for ambit processes, which will not be demonstrated in this work,
are investigated in \cite{BEV}.

The paper is structured as follows. We introduce the main mathematical setting and discuss some
basic probabilistic properties of ambit fields in Section \ref{sec2}. Section \ref{sec3} 
gives a review on integration concepts with respect to ambit fields.  Section \ref{sec5} 
is devoted to limit theorems for high frequency statistics of ambit fields and their statistical 
applications.

\section{The model and basic probabilistic properties} \label{sec2}
\setcounter{equation}{0}
\renewcommand{\theequation}{\thesection.\arabic{equation}}
Throughout this paper all stochastic processes are defined on a given probability space $(\Omega, \mathcal F, \mathbb P)$.
We first consider a subclass of an ambit field defined via
\begin{align} \label{statambfield}
  X_{t}(x)=\mu +\int_{A+(t,x)}g(t-s, x-\xi )\sigma_{s}(\xi )L(\mathrm{d}s\mathrm{d}
  \xi )+\int_{D+(t,x)}q(t-s,x-\xi )a_{s}(\xi )\mathrm{d}s\mathrm{d}\xi,
\end{align}
where $A,D\subset \R \times \R^d$ are fixed ambit sets. This specific framework allows for modelling stationary
random fields, namely when the field $(\sigma , a)$ is stationary and independent of the driving L\'evy basis $L$,
$X$ turns out to be stationary too. Another potentially important features of a random field, such as e.g. isotropy
or skewness,
can be easily modelled via an appropriate choice of the deterministic kernels $g$ and $q$, and the stochastic field $(\sigma , a)$.

Specializing even further, we introduce a \textit{L\'evy semi-stationary process} given by
\begin{align} \label{lss}
  Y_{t}=\mu +\int_{-\infty}^t g(t-s)\sigma_{s} L(\mathrm{d}s)+
\int_{-\infty}^t q(t-s)a_{s} \mathrm{d}s,
\end{align}   
where now $L$ is a two-sided one dimensional L\'evy motion and $A=D=(-\infty ,0)$, which is a purely temporal analogue
of the ambit field at \eqref{statambfield}. It is indeed this special subclass of ambit fields, which will be the focus
of our interest, as many mathematical problems have to be solved for L\'evy semi-stationary processes first before
transferring the principles to more general random fields. In the following we will check the validity 
and discuss some basic properties of ambit fields. Note that we sometimes use $dL_s$ instead of $L(\mathrm{d}s)$
in the setting of \eqref{lss}. We also remark that the meaning of some notations may change from section to section.

\subsection{Definition of the stochastic integral}
The first natural question when introducing an ambit field is the definition of the first integral
appearing in the decomposition \eqref{ambfield}. Here we recall two classical approaches of Rajput and Rosinski \cite{RR}
and Walsh \cite{W}.

\subsubsection{Rajput and Rosinski theory}
The first simplified mathematical problem is the definition of stochastic integral in case of a \textit{deterministic}
intermittency field $\sigma$. In this situation we may study a more general integral
\[
\int_A f d L,
\]
where $L$ is a L\'evy basis on a $\delta$-ring 
$\mathcal {S}$, $f: (S, \sigma (\mathcal {S}))\rightarrow (\R , \mathcal{B}(\R))$
a measurable real valued function and $A\in \sigma (\mathcal {S})$, which is the main object of a seminal paper \cite{RR}.
We will briefly recall the most important results of this work. By definition of a L\'evy basis the law of $(L(A_1), \ldots,
L(A_d))$, $A_i\in \mathcal {S}$, is infinitely divisible, the random variables 
$L(A_1), \ldots, L(A_d)$ are independent when the sets $A_1, \ldots, A_d$ are disjoint, and
\[
L\left(\cup_{i=1}^\infty A_i \right) = \sum_{i=1}^\infty L(A_i) \qquad \mathbb{P}-\text{almost surely}
\] 
for disjoint $A_i$'s with $\cup_{i=1}^\infty A_i\in \mathcal {S}$. Recalling the characteristic triplet $(v_1,v_2,v_3)$
from the L\'evy-Khinchin representation \eqref{charlevy}, the \textit{control measure} $\lambda$ of $L$ is defined via
\[
\lambda (A): = |v_1|(A) + v_2(A) + \int_{\R} \min (1, x^2) v_3(dx, A), \qquad A\in \mathcal{S},
\]  
where $|v_1|$ denotes the total variation measure associated with $v_1$. In this subsection we will use
the truncation function
\[
\tau (z):= \left\{ 
\begin{array} {cc}
z: & \|z\|\leq 1 \\
z/\|z\|: & \|z\|>1
\end{array}
\right.
\]
Now, for any \textit{simple function} 
\[
f(x)=\sum_{i=1}^d a_i 1_{A_i} (x), \qquad a_i\in \R, A_i\in \mathcal {S},
\]
the stochastic integral is defined as
\[
\int_{A} f dL:= \sum_{i=1}^d a_i L(A\cap A_i), \qquad A\in \sigma (\mathcal{S}).
\] 
The extension of this definition is as follows. 

\begin{defi} \label{def1} \rm
A measurable function $f: (\mathcal {S}, \sigma (\mathcal {S}))\rightarrow (\R , \mathcal{B}(\R))$ is called $L$-\textit{integrable} 
if there exists a sequence of simple functions $(f_n)_{n\geq 1}$ such that
\begin{itemize}
\item[(i)] $f_n \rightarrow f$ $\lambda$-almost surely. 
\item[(ii)] For any $A\in \sigma (\mathcal{S})$ the sequence $(\int_A f_n dL)$ converges in probability.
\end{itemize} 
In this case the stochastic integral is defined by
\[
\int_{A} f dL:= \mathbb{P}-\lim_{n\rightarrow \infty } \int_A f_n dL. 
\]
\end{defi}
Although this definition is quite intuitive, it does not specify the class of $L$-integrable functions explicitly. The next theorem,
which is one of the main results of \cite{RR}, gives an explicit condition on the $L$-integrability of a function $f$. 

\begin{theo} \label{th1}
(\cite[Theorem 2.7]{RR}) Let $f: (S, \sigma (\mathcal {S}))\rightarrow (\R , \mathcal{B}(\R))$ be a measurable function. Then $f$ is $L$-integrable
if and only if the following conditions hold:
\[
\int_{S} U(f(s),s) \lambda(ds)<\infty, \quad \int_{S} f^2(s)v_2(s) \lambda(ds)<\infty, \quad
\int_{S} V_0(f(s),s) \lambda(ds)<\infty,   
\]
where 
\begin{align*}
U(u,s) &:= u v_1(s) + \int_{\R} ( \tau (xu) - u\tau(x) ) v_3(dx,s), \\
V_0(u,s)&:= \int_{\R} \min (1, |xu|^2) v_3(dx,s). 
\end{align*}
Furthermore, the real valued random variable $X=\int_{ S} f dL$ is infinitely divisible with L\'evy-Khinchin representation
\[
\log \E[\exp(iu X)]= iu v_1(f) - \frac 12 u^2 v_2(f) + \int_{\R} \left( \exp(iuy) - 1- iuy 1_{[-1,1]}(y) \right) v_3^f(dy),
\]  
where
\begin{align*}
v_1(f)&= \int_{S} U(f(s),s) \lambda (ds), \\
v_2(f)&= \int_{S} f^2(s)v_2(s) \lambda(ds), \\
v_3^f(B)&= v_3 \left\{(x,s)\in \R\times S:~ xf(s)\in B\setminus \{0\}\right\}, \qquad B\in \mathcal{B}(\R).
\end{align*} 
\end{theo}

\begin{rem} \rm \label{rem1}
The most basic example of a null spatial ambit field is the L\'evy moving average process given by
\begin{align} \label{lma}
X_t=\int_{-\infty}^t g(t-s) dL_s,
\end{align}
where $L$ is a two-sided L\'evy motion with $\E L_t=0$. In this situation the control measure $\lambda$ is just the Lebesgue measure
and the sufficient conditions from Theorem \ref{th1} translate to 
\[
 \int_0^\infty g^2(s) ds <\infty, \quad \int_0^{\infty } \left( 
\int_{\R} \min(|x g(s)|, |x g(s)|^2) v_3(dx)\right) ds <\infty,  
\]
where $v_3$ is the L\'evy measure of $L$. \qed
\end{rem}
Theorem \ref{th1} gives a precise condition for existence of the stochastic integral. However, it does not say anything about the
existence of moments. To study this question let us assume for the rest of this subsection that
\[
\E[|L(A)|^q]<\infty \qquad \forall A\in \mathcal{S},
\]  
for some $q>0$. Clearly, in general we can not expect the existence of moments higher than $q$ for the integral $\int_{S}
f dL$, but we may study the existence of $p$th moment with $p\leq q$. For this purpose we introduce the function $\Phi_p:
\R \times S \rightarrow \R$
via
\[
\Phi_p(u,s):= U^{\star} (u,s) + u^2 v_2(s) + V_p(u,s),
\] 
where 
\begin{align*}
U^{\star} (u,s)&:= \sup_{c\in [-1,1]} |U(cu,s)|, \\
V_p(u,s)&:= \int_{\R} \left(|ux|^p 1_{\{|ux|>1\}} + |ux|^2 1_{\{|ux|\leq 1\}} \right) v_3 (dx,s),
\end{align*}
and the function $U$ is introduced in Theorem \ref{th1}. Then 
\[
\| f\|_{\Phi_p}:= \inf \left \{c>0:~ \int_{S} \Phi_p( c^{-1} |f(s)|,s) \lambda (ds)  \right \}
\]
defines a norm and the vector space 
\[
\mathbb L^{\Phi_p}:= \left\{f:~  \int_{S} \Phi_p(  |f(s)|,s) \lambda (ds)  <\infty \right\}
\]
equipped with $\| \cdot \|_{\Phi_p}$ is the so called \textit{Musielak-Orlicz space} (see \cite{RR} for details).
The next result gives a connection between the existence of $p$th norm of the stochastic integral and the finiteness
of  $\| f\|_{\Phi_p}$.

\begin{theo} \label{th2}
(\cite[Theorem 3.3]{RR}) Let $p\in [0,q]$. Then it holds that 
\[
\E\left[ \left|\int_{S}
f dL \right|^p \right]<\infty \Longleftrightarrow \| f\|_{\Phi_p}<\infty. 
\]
\end{theo}

\begin{rem} \rm \label{rem2}
Let $X$ be a L\'evy moving average process as defined at \eqref{lma}, where the driving L\'evy motion has characteristic 
triplet $(0,0, \text{const.} |x|^{-1-\beta }dx)$, i.e. $L$ is a symmetric $\beta$-stable process ($\beta \in (0,2)$)
without drift. In this case the $p$th absolute moment of $X_t$ exists whenever $g\in L^{\beta} (\R_{\geq 0})$ and 
$p<\beta$. Moreover it holds that
\[
C_1 \E \left[ \left| \int_0^\infty |g(x)|^{\beta} dx \right| \right]^{1/\beta } \leq \E[|X_t|^p]^{1/p} \leq
C_2 \E \left[ \left| \int_0^\infty |g(x)|^{\beta} dx \right| \right]^{1/\beta } 
\]
for some positive constants $C_1, C_2$. \qed
\end{rem}

\subsubsection{Walsh approach}
The integration concept proposed by Walsh \cite{W} is, in some sense, an extension of It\^o's theory
to spatio-temporal setting. As in the classical It\^o calculus the martingale theory and isometry play 
an essential role.  

Here we review the main ideas of Walsh approach. 
In the following we will concentrate on integration with respect to a L\'evy basis on \textit{bounded} domains (the theory can be extended 
to unbounded domains in a straightforward manner, see \cite{W}). 
Let $\mathcal B_b(\R^d)$ the  Borel $\sigma$-field generated by bounded Borel sets on $\R^d$ and 
$L$  a L\'evy basis on $[0,T] \times S$ with $S\in \mathcal B_b(\R^d)$.
For $S\supseteq A\in \mathcal B_b(\R^d)$ we introduce the notation 
\[
L_t(A):= L((0,t]\times A).
\] 
Unlike in the approach of Rajput and Rosinski, which does not require any moment assumptions 
(which however deals only with deterministic integrands), the basic
assumption of Walsh concept is: 
\begin{center}
For any $A\in \mathcal B_b(\R^d)$ it holds that $\E[L_t(A)]=0$ and $L_t(A)\in L^2(\Omega, \mathcal F, \mathbb P)$.
\end{center}
Now, we introduce a right continuous filtration $(\mathcal F_t)_{t\geq 0}$ via
\[
\mathcal F_t := \cap_{\epsilon >0} \mathcal F^0_{t+\epsilon} \quad \text{with} \quad 
F^0_{t}= \sigma \{ L_s(A):~ S\supseteq A\in \mathcal B_b(\R^d), 0<s\leq t\}  \vee \mathcal N,
\]
where $\mathcal N$ denotes the $\mathbb{P}$-null sets of $\mathcal F$. Since $L$ is a L\'evy basis, the stochastic
field $(L_t(A))_{t\geq 0, S\supseteq A\in \mathcal B_b(\R^d)}$ is an \textit{orthogonal martingale measure} with respect
to $(\mathcal F_t)_{t\geq 0}$ in the sense of Walsh,
i.e. $(L_t(A))_{t\geq 0}$ is a square integrable  martingale with respect to the filtration 
$(\mathcal F_t)_{t\geq 0}$ for any $A\in \mathcal B_b(\R^d)$ and the random variables $L_t(A), L_t(B)$ are independent
for disjoint sets $A,B\in \mathcal B_b(\R^d)$.

In the next step we introduce the \textit{covariance measure} $Q$ via
\[
Q([0,t]\times A):= \langle L(A) \rangle_t, \qquad A\in \mathcal B_b(\R^d), 
\]
and the associated $L^2$-norm
\[
\|\psi \|_{Q}^2 := \E \left[ \int_{[0,T] \times S} \psi^2 (t, \xi) Q(dt,d \xi)  \right]
\]
for a random field $\psi$ on $[0,T] \times S$. Now, we are following the It\^o's program: We call $\psi$ an
\textit{elementary random field} if it has the form
\[
\psi(\omega , s, \xi)=  X(\omega ) 1_{(a,b]} (s) 1_A(\xi),
\]
where $X$ is a bounded $\mathcal F_a$-measurable random variable and $A\in \mathcal B_b(\R^d)$. For such an elementary 
random field the integral with respect to a L\'evy basis is defined by ($t\leq T$)
\[
\int_0^t \int_B \psi( s, \xi) L(ds, d \xi):= X (L_{t\wedge b} (A\cap B) - L_{t\wedge a} (A\cap B)), \qquad B\in \mathcal B_b(\R^d),~
A,B\subseteq S.
\]
By linearity this definition can be extended to the linear span given by elementary random fields. The $\sigma$-field $\mathcal P$
generated by elementary random fields is called \textit{predictable}, and the space $\mathcal P_L= L^2(\Omega \times [0,T] \times
S, \mathcal P, Q )$ equipped with $\| \cdot \|_Q$ is a Hilbert space. Furthermore, the space of elementary 
functions is dense in $\mathcal P_L$. Thus, for any random field $\psi \in \mathcal P_L$, we may find a sequence 
of elementary random fields $\psi_n$ with $\| \psi_n-\psi  \|_Q\rightarrow 0$ and 
\[
\int_0^t \int_B \psi( t, \xi) L(dt, d \xi):= \lim_{n\rightarrow \infty} \int_0^t \int_B \psi_n( t, \xi) L(dt ,d \xi)
\qquad \text{in } L^2(\Omega , \mathcal F, \mathbb P).
\] 
Moreover, the It\^o type isometry 
\[
\E \left[ \left| \int_0^T \int_S \psi( s, \xi) L(ds, d \xi) \right|^2 \right] = \|\psi \|_{Q}^2
\]
is satisfied by construction.

\begin{rem} \label{rem3} \rm
The Walsh approach can be used to define ambit fields in \eqref{ambfield} only when the driving L\'evy basis 
is square integrable, which excludes e.g. $\beta$-stable processes ($\beta \in (0,2)$). The recent work 
of \cite{CK} combines the ideas of Walsh \cite{W} and Rajput and Rosinski \cite{RR} to propose an integration 
concept for random integrands and general L\'evy bases. It relies on an earlier work \cite{BJ} by Bichteler and Jacod. 
A general reference for comparison of various integration concepts is \cite{DQ}. \qed  
\end{rem}

\subsection{Is an ambit field a semimartingale?}
Semimartingales is nowadays a well studied object in probability theory. From theoretical perspective it is important to understand
whether a given null spatial subclass of ambit fields is a semimartingale or not, since in this case one may better study its fine
structure properties. Furthermore, limit theory for high frequency statistics, which will be the focus of our discussion in Section 
\ref{sec5}, has been investigated in great generality in the framework of semimartingales; we refer to \cite{JP} for a comprehensive
study. Thus, new asymptotic results are only needed for ambit fields, which do not belong to the class of semimartingales.

For simplicity of exposition we concentrate on the study of \textit{Volterra type} processes
\begin{align} \label{volterra}
X_t = \int_{-\infty}^t g(t,s) dL_s,
\end{align}  
where $L$ is a L\'evy motion with $L_0=0$, which is obviously a purely temporal subclass of ambit fields. The semimartingale property
can be studied with respect to the filtration generated by the driving L\'evy motion $L$ or with respect to the natural filtration
(which is a harder task). For any $t\geq 0$, we define
\[
\mathcal F^L_t :=\sigma ( L_s:~s\leq t), \qquad \mathcal F^X_t :=\sigma ( X_s:~s\leq t).
\]  
In the following we will review the main studies of the semimartingale property.

\subsubsection{Semimartingale property with respect to $(\mathcal F^L_t)$}
Semimartingale property of various subclasses of Volterra type processes $X$ have been studied in the literature, see e.g.
\cite{BP,BLS,K} among many others. However, in this subsection we closely follow the recent results of \cite[Section 4]{BR}. We note
that the original work \cite{BR} contains a study of more general processes, which are specialized in this paper to models
of the form \eqref{volterra}.

Let $L$ be a L\'evy process with characteristic triplet $(a,b^2, \nu )$. In the following we consider \textit{stationary
increments  moving average} model of the type
\begin{align} \label{SIMMA}
X_t= \int_{\R} [g_1(t-s) -g_0(-s)] dL_s, 
\end{align}   
where $g_0,g_1$ are measurable functions satisfying $g_0(x)=g_1(x)=0$ for $x<0$. This obviously constitutes a subclass
of \eqref{volterra}. Integrability conditions can be directly extracted from Theorem \ref{th1}. Notice that by construction 
the process $X$ has stationary increments, which explains the aforementioned notion. The main result of this subsection is 
\cite[Theorem 4.2]{BR}.

\begin{theo} \label{th3.0}
Assume that the process $X$ is defined as in \eqref{SIMMA} and the following conditions are satisfied:
\begin{align*}
&\int_0^{\infty} |g'_1 (s)|^2 ds <\infty \qquad(\text{when } b^2>0), \\
& \int_0^{\infty} \int_{\R} \min(|xg'_1 (s)|, |xg'_1 (s)|^2) \nu(dx) ds<\infty.    
\end{align*}
Then $X$ is a semimartingale with respect to the filtration $(\mathcal F^L_t)$.
\end{theo}
Theorem \ref{th3.0} gives sufficient conditions for the semimartingale property of the process $X$. In certain
special cases these conditions are also necessary as the following result from \cite[Corollary 4.8]{BR} shows.

\begin{theo} \label{th3}
Assume that the process $X$ is defined as in \eqref{SIMMA}, $L$ has infinite variation on compact sets 
and $L_1$ is either square integrable 
or has a regular varying distribution at $\infty$ with index $\beta \in [-2,-1)$. Then $X$ is a semimartingale with respect to the filtration $(\mathcal F^L_t)$ if and only if the following conditions are satisfied:
\begin{align*}
&\int_0^{\infty} |g'_1 (s)|^2 ds <\infty \qquad(\text{when } b^2>0), \\
& \int_0^{\infty} \int_{\R} \min(|xg'_1 (s)|, |xg'_1 (s)|^2) \nu(dx) ds<\infty.    
\end{align*}
In this case it has the decomposition 
\[
X_t = X_0+ g_1(0) L_t + \int_0^t \left( \int_{\R} g'_1(s-u) dL_u \right) ds. 
\]  
\end{theo}  
When the driving process $L$ is a symmetric $\beta$-stable L\'evy motion with $\beta \in (1,2 )$, i.e. the characteristic triplet
is given via $(0,0, \text{const.} |x|^{-1-\beta }dx)$, the condition of Theorem \ref{th3} translates to 
\[
\int_0^{\infty} |g'_1 (s)|^{\beta} ds <\infty.
\]

\begin{rem} \rm \label{rem4}
Another important subclass of stationary increments  moving average models, which will be the object of investigation in Section
\ref{sec5}, is a \textit{fractional L\'evy motion}. A fractional L\'evy motion is defined as
\begin{align} \label{flm}
X_t= \int_{-\infty}^t [(t-s)^{\alpha}_+ - (-s)^{\alpha}_+] dL_s.
\end{align} 
When $L$ is a Brownian motion, this is one of the representation of a fractional Brownian motion
with Hurst parameter $H=\alpha +1/2$ (with $\alpha \in (-1/2,1/2)$).
Assume now that $\alpha >0$. Then $X$ is a semimartingale with respect to the filtration $(\mathcal F^L_t)$ if and only if $b^2=0$,
$\alpha  \in (0, 1/2)$ and 
\[
\int_{\R} |x|^{1/(1-\alpha )} \nu (dx)<\infty. 
\]  
We refer to \cite[Proposition 4.6]{BR} for details. \qed 
\end{rem}

\begin{rem} \rm
Sufficient conditions for the semimartingale property for general L\'evy 
semi-stationary processes defined at \eqref{lss} can be deduced from 
a seminal work \cite{Pr} on Volterra type equations. \qed  
\end{rem}

\begin{rem} \rm \label{rem5}
In physical applications it may appear that an ambit field is observed along a curve in time-space. Let us consider
for simplicity an ambit field of the form
\[
X_{t}(x)= \int_{A_{t}(x)}g(t,s,x,\xi ) L(\mathrm{d}s,\mathrm{d} \xi ),
\] 
where $L$ is a L\'evy basis on $\R_{\geq 0} \times  \R^d$. Let $\theta =(\theta_1, \theta_2) 
: [0,T]\rightarrow \R_{\geq 0} \times  \R^d$ be a curve in time-space. Then the observed process is given by
\[
Y_t = X_{\theta_1 (t)}(\theta_2(t)).
\]
The semimartingale property of the process $Y$ is still an open problem. \qed 
\end{rem}

\subsubsection{Semimartingale property with respect to $(\mathcal F^X_t)$} 
Showing the semimartingale property with respect to the natural filtration $(\mathcal F^X_t)$ is by far a more delicate issue.
In this subsection we restrict our attention to Volterra processes of the form
\[
\int_{\R} [g_1(t-s) -g_0(-s)] dW_s.
\]
where $W$ is a Brownian motion, since, to the best of our knowledge, little is known for general driving L\'evy processes.
Here $g_0,g_1$ are measurable functions satisfying the integrability condition $s\mapsto g_1(t-s) -g_0(-s)\in L^2(\R)$
for any $t\in \R$ (this insures the existence of the integral). In the case $g_0=0$ the paper \cite{JY} provides
necessary and sufficient conditions for the semimartingale property with respect to the natural filtration $(\mathcal F^X_t)$.
The work \cite{B} extends the results to include the model defined above. The methodology of proofs is based upon 
Fourier transforms and Hardy functions.

We need to introduce some notation. For a function $h:\R \rightarrow \R$, $\widehat{h}:\R\rightarrow \mathbb C$ denotes
its Fourier transform, i.e.
\[
\widehat{h}(t): =\int_{\R} \exp(itx) h(x)dx.
\]
For a function $f: \R\rightarrow S^1$, where $S^1$ denotes the unit circle on a complex plane, satisfying $\overline f(\cdot)
= f(-\cdot)$ we define the function $\widetilde f:\R \rightarrow \R$ via
\[
\widetilde f (t):= \lim_{a\rightarrow \infty} \int_{-a}^a \frac{\exp(its) - 1_{[-1,1]}(s)}{is} f(s) ds.
\]
The main result of this subsection is \cite[Theorem 3.2]{B}.

\begin{theo} \label{th4}
The Gaussian process $X$ is a semimartingale with respect to the natural filtration $(\mathcal F^X_t)$ if and only if 
the following conditions are satisfied:
\begin{itemize}
\item[(i)] The function $g_1$ has the representation 
\[
g_1(t) = b+a \widetilde{f}(t) + \int_0^t \widehat{f \hat{h}}(s) ds \qquad \text{for Lebesgue almost all }t\in \R,
\]
where $a,b\in \R$, $f: \R\rightarrow S^1$ is a measurable function with $\overline f(\cdot)
= f(-\cdot)$, and $h\in L^2(\R)$ is $0$ whenever $a\not = 0$.
\item[(ii)] Set $\zeta = \widehat{f (\widehat{g_1-g_0})}$. When $a\not= 0$ then
\[
\int_0^r \left( \frac{|\zeta(s)|}{\sqrt{\int_s^{\infty } \zeta^2(u) du}} \right) ds<\infty \qquad \forall r>0,
\] 
where $0/0:=0$.
\end{itemize} 
\end{theo}   
The second part of \cite[Theorem 3.2]{B} gives the canonical decomposition of $X$ in case, where the above conditions (i)
and (ii) are satisfied. If $X_0=0$ one may choose $a,b,h$ and $f$ such that the canonical decomposition $X_t= M_t+A_t$
is given via
\[
M_t= a \int_{\R} (\widetilde{f}(t-s) - \widetilde{f}(-s)) dW_s, \qquad A_t= \int_0^t \left( \int_{\R}
 \widehat{f \hat{h}}(s-u) dW_u \right) ds.
\]  
In this decomposition the martingale part $M$ is a Wiener process with scaling parameter $(2\pi a)^2$. Finally, let us remark
that in the special case $g_1=g_0$, as e.g. in the fractional Brownian motion setting, the condition (ii) of Theorem \ref{th4}
is trivially satisfied.

\subsection{Fine properties of L\'evy semi-stationary processes} \label{sec2.3}
In Section \ref{sec5} we will study the limit theory for \textit{power variation} of L\'evy semi-stationary processes
\[
\sum_{i=1}^{[t/\Delta_n]} |Y_{i\Delta_n} - Y_{(i-1)\Delta_n}|^p \qquad \text{with } \Delta_n\rightarrow 0, 
\] 
where the process $Y$ is defined at \eqref{lss}. A first step towards understanding the  first and second order asymptotics for
this class of statistics is to determine the fine scale behaviour of L\'evy semi-stationary processes. For simplicity of exposition
we start our discussion with L\'evy moving average processes defined at \eqref{lma}, i.e. 
\[
X_t=\int_{-\infty}^t g(t-s) dL_s,
\] 
which constitute a subclass of L\'evy semi-stationary processes with constant intermittency and zero drift. We restrict ourselves
to symmetric $\beta$-stable  L\'evy processes $L$ with $\beta \in (0,2]$ and zero drift (thus, the Brownian motion is included). 

The most interesting class of weight functions $g$ in physical applications is given via
\begin{align} \label{basicex}
g(x)=x^{\alpha} f(x) \qquad \text{for } x>0,
\end{align}
and $g(x)=0$ for $x\leq 0$, where $f: \R_{\geq 0} \rightarrow \R$ is a smooth function with $f(0) \not=0$ 
decaying fast enough at infinity 
to ensure the existence of the integral (see Remark \ref{rem1}). When $\beta \in (0,2)$, i.e. $L$ is a pure jump process, 
we further always assume  that $\alpha >0$, since $\alpha <0$ leads to explosive behaviour of $X$ near jump times of $L$. Hence,
for any $\beta \in (0,2]$, the L\'evy moving average process $X$ is continuous since $g(0)=0$ when $\beta <2$. A formal
differentiation leads to the identity
\[
dX_t = g(0+) dL_t +  \left( \int_{-\infty}^t g'(t-u) dL_u \right) dt.
\]  
According to Theorem \ref{th3} this identity indeed makes sense when (a) $\beta =2$, $g(0+)<\infty$ and $g'\in L^2(\R_{\geq 0})$,
or (b) $\beta \in (0,2)$ and
\[
\int_{0}^\infty |g'(x)|^{\beta} dx<\infty. 
\]
In this case the L\'evy moving average process $X$ is an It\^o semimartingale and the law of large numbers for its 
power variation is well understood (see e.g. the monograph \cite{JP}). In Section \ref{sec5} we will be specifically interested
in situations where $X$ is not a semimartingale. It is particularly the case under the conditions 
\begin{align*}
& \beta=2 \quad \text{and} \quad \alpha \in (-1/2, 1/2) \quad \text{or} \\
& \beta \in (0,2)  \quad \text{and} \quad \alpha \in (0, 1-1/\beta),
\end{align*}  
since the above integrability condition for the function $g'$ is not satisfied near $0$. Under these conditions the derivative
of the function $g$ explodes at $0$. For a small $\Delta>0$, we intuitively deduce the following approximation
for the increments of $X$:

\begin{align*}
X_{t+\Delta} - X_t &= \int_{\R} [g(t+\Delta -s) - g(t -s)] dL_s \\
& \approx \int_{t+\Delta -\epsilon }^{t+\Delta} [g(t+\Delta -s) - g(t -s)] dL_s \\
& \approx f(0) \int_{t+\Delta -\epsilon }^{t+\Delta} [(t+\Delta -s)_{+}^\alpha  - (t -s)_{+}^\alpha ] dL_s \\
& \approx f(0) \int_{\R} [(t+\Delta -s)_{+}^\alpha  - (t -s)_{+}^\alpha ] dL_s = \widetilde{X}_{t+\Delta} - \widetilde{X}_t,
\end{align*} 
where 
\begin{align} \label{flm1}
\widetilde{X}_t = f(0) \int_{\R} [(t-s)_{+}^\alpha  - ( -s)_{+}^\alpha ] dL_s,
\end{align} 
and $\epsilon >0$ is an arbitrary small real number with $\epsilon \gg \Delta$. The formal proof of this first order approximation
relies on the fact that the weight $g(t+\Delta -s) - g(t -s)$ attains asymptotically highest values when $s\approx t$, since
$g'$ explodes at $0$. The stochastic process $\widetilde{X}$ is called a \textit{fractional $\beta$-stable L\'evy motion} and its
properties have been studied in several papers, see e.g. \cite{BCI,M} among others 
(for $\beta =2$ it is just an ordinary scaled fractional Brownian
motion with Hurst parameter $H=\alpha +1/2$). 
In particular, $\widetilde{X}$ has stationary increments, it is 
$(\alpha +1/\beta)$-self similar with symmetric $\beta$-stable marginals. 

The key fact to learn from this approximation is that, under above assumptions on $g$,  
the fine structure of a L\'evy moving average process $X$ with symmetric $\beta$-stable driver $L$ is similar to the fine structure
of a fractional $\beta$-stable L\'evy motion $\widetilde{X}$. Thus, under certain conditions, one may transfer the asymptotic 
theory for power variation of $\widetilde{X}$ to the corresponding results for power variation of $X$. The limit theory
for power variation of $\widetilde{X}$ is sometimes easier to handle than the original statistic due to stationarity of increments
of $\widetilde{X}$ and their self similarity property, which allows to transform the original triangular observation scheme
into a usual one when studying distributional properties (for the latter, ergodic limit theory might apply). Indeed, it is one 
method of proofs of laws of large numbers presented in Theorem \ref{th7}(ii) of Section \ref{sec5}.

\begin{rem} \label{rem6} \rm
As mentioned above the fractional L\'evy motion defined at \eqref{flm1} is, up to a scaling factor, a fractional 
Brownian motion $B^H$ with Hurst parameter $H=\alpha +1/2 \in (0,1)$ when $\beta =2$. Due to self similarity property 
of $B^H$ it is sufficient to study the asymptotic behaviour of the statistic 
\[
\sum_{i=1}^{[t/\Delta_n]} |B_i^H - B_{i-1}^H|^p, \qquad p>0,
\]  
to investigate the limit theory for power variation of $B^H$. 
Below we review some classical results of \cite{BM83,T}. First of all, we obtain the convergence
\begin{align} \label{mp}
\Delta_n \sum_{i=1}^{[t/\Delta_n]} |B_i^H - B_{i-1}^H|^p \ucp m_p t, \qquad m_p:= \E[|\mathcal N(0,1)|^p],
\end{align}
where $Z^n \ucp Z$ stands for uniform convergence in probability on compact intervals, i.e. $\sup_{t\in [0,T]} |Z^n_t-Z_t| \toop 0$. 
The associated weak limit theory depends on the correlation kernel of the fractional Brownian noise 
and the \textit{Hermite rank} of the function $h(x)=|x|^p - m_p$. Recall that the correlation kernel of the fractional Brownian noise
is given via
\[
\rho (j):=\text{corr} (B_1^H-B_0^H, B_{j+1}^H-B_j^H)= \frac{1}{2} \Big(|j+1|^{2H} - 2|j|^{2H} + |j-1|^{2H} \Big).
\]
The Hermite expansion of the function $h$ is defined as
\[
h(x)=|x|^p - m_p = \sum_{l=2}^{\infty} a_l H_l(x),
\]
where $(H_l)_{l\geq 0}$ are Hermite polynomials, i.e. 
\[
H_0=1 \quad \text{and} \quad H_l = (-1)^l \exp(x^2/2) \frac{d}{dx^l} \{-\exp(x^2/2)\} \quad \text{for } l\geq 1.
\]
The Hermite rank of $h$ is the smallest index $l$ with $a_l \not =0$, which is $2$ in our case. The condition for the validity
of a central limit theorem associated with \eqref{mp} is then
\[
\sum_{j=1}^\infty \rho^2 (j)<\infty, 
\] 
where the power $2$ indicates the Hermite rank of $h$. This condition holds if and only if $H\in (0, 3/4)$. For $H>3/4$ 
the limiting process is non-central. More precisely, the following functional limit theorems hold:
\begin{align*}
0<H<3/4 &: \qquad \Delta_n^{-1/2}\left( 
\Delta_n \sum_{i=1}^{[t/\Delta_n]} |B_i^H - B_{i-1}^H|^p - m_p t \right) \schw v_p W'_t, \\[1.5 ex]
H=3/4 &:   \qquad (\Delta_n \log{\Delta_n^{-1}})^{-1/2}\left( 
\Delta_n \sum_{i=1}^{[t/\Delta_n]} |B_i^H - B_{i-1}^H|^p - m_p t \right) \schw \tilde{v}_p W'_t, 
\\[1.5 ex]
3/4<H<1 &:   \qquad \Delta_n^{2H-2}\left( \Delta_n \sum_{i=1}^{[t/\Delta_n]} |B_i^H - B_{i-1}^H|^p - m_p t \right) \schw Z_t,
\end{align*}
where the weak convergence takes place on $\mathbb D([0,T])$ equipped with the uniform topology, $W'$ denotes a Brownian motion
and $Z$ is a Rosenblatt process (see e.g. \cite{T}). Finally, the constants $v_p$ and $\tilde{v}_p$ are given by
\begin{align*}
v_p &:= \sum_{l=2}^\infty l! a_l^2 \Big(1+ 2 \sum_{j=1}^\infty \rho^l(j) \Big), \\
\tilde{v}_p  &:= 2 \lim_{n\rightarrow \infty} \frac{1}{\log n} \sum_{j=1}^{n-1} \frac{n-k}{n} \rho^2(j) \cdot
\sum_{l=2}^\infty l! a_l^2.
\end{align*}
\qed
\end{rem}

\begin{rem} \label{rem7} \rm
When $\beta \in (0,2)$ and $\alpha \in (0, 1-1/\beta)$ we deduce via the $(\alpha +1/\beta )$-self similarity 
property and the strong ergodicity of the fractional $\beta$-stable L\'evy process $\widetilde{X}$
\[
\Delta_n^{1-p(\alpha +1/\beta )}\sum_{i=1}^{[t/\Delta_n]} |\widetilde{X}_{i\Delta_n} - \widetilde{X}_{(i-1)\Delta_n}|^p
\ucp c_p t, \qquad c_p:= \E[|\widetilde{X}_1- \widetilde{X}_0|^p] 
\] 
whenever $p<\beta$. For $p>\beta $ the constant $c_p$ is infinite and non-ergodic limits appear, see Theorem \ref{th7} in Section \ref{sec5}. \qed
\end{rem}

\begin{rem} \label{rem8} \rm
Once we have proved the law of large numbers for power variation of the basic process $\widetilde{X}$ as in Remarks
\ref{rem6} and \ref{rem7} (or for the L\'evy moving average process $X$ defined at \eqref{lma}), the main principles
of the proof usually transfer to the integral process
\[
I_t=\int_0^t \sigma_s d \widetilde{X}_t,
\] 
whenever the latter is well defined. As it was shown in e.g. \cite{BCP09,CNW} the Bernstein's blocking technique can be applied
to deduce the law of large numbers for power variation of the process $I$. For instance, when $\widetilde{X}=B^H$
is a fractional Brownian motion with Hurst parameter $H$, it holds that 
\[
\Delta_n^{1-pH}\sum_{i=1}^{[t/\Delta_n]} |I_{i\Delta_n} - I_{(i-1)\Delta_n}|^p \ucp m_p \int_0^t |\sigma_s|^p ds,
\] 
where $\sigma$ is a stochastic process, which has finite $q$-variation with $q<1/(1-H)$ (see \cite{CNW}). Quite often
the same asymptotic result holds also for a subclass of ambit processes given by
\[
Y_t=\int_{-\infty}^t g(t-s) \sigma_s dW_s,
\]
where $W$ is a Brownian motion (cf. Section \ref{sec5}). The reason is again a similar fine structure of the processes
$I$ and $Y$ (see e.g. \cite[Section 2.2.3]{CHPP13} for a detailed exposition). 

Transferring a central limit theorem for the power variation of  the driver $B^H$ to the integral process $I$
(and also in case of the process $Y$) is a more delicate 
issue. Apart from further assumptions on the integrand $\sigma$ a more technical and 
precise treatment of the Bernstein's blocking technique
is required. We refer to a recent work \cite{CNP14} for a detailed description of such a method, which relies on fractional calculus.
\qed      
\end{rem}

\section{Integration with respect to ambit fields} \label{sec3}
\setcounter{equation}{0}
\renewcommand{\theequation}{\thesection.\arabic{equation}}
In this section we will discuss the integration concepts with respect to ambit processes of the type
\[
X_t= \int_0^t g(t,s)\sigma_s dL_s,
\] 
where $\sigma$ is a stochastic intermittency process and $L$ is a L\'evy motion. Our presentation is mainly based upon the recent work
\cite{BBPV}, where Malliavin calculus is applied to define the stochastic integral. The introduction of the integral

\begin{align} \label{integral}
\int_0^t Z_s dX_s,
\end{align} 
for a stochastic integrand $Z$, strongly depends on whether the driving process $L$ is a Brownian motion or a pure jump L\'evy motion, 
since the main notions of Malliavin calculus differ in those two cases. An alternative way of defining the stochastic integral
at \eqref{integral} without imposing $L^2$-structure of the integrand $Z$ is proposed in \cite{BBS}. The authors apply white noise analysis
to construct the integral in the situation, where $X$ is driven by a Brownian motion.

\subsection{Integration with respect to ambit processes driven by Brownian motion}
Before we present the definition of the integral at \eqref{integral} for $L=W$, we start by introducing the main notions of Malliavin calculus on Gaussian spaces.
The reader is referred to the monograph \cite{N} for any unexplained definition or result. 

Let $\HH$ be a real separable Hilbert space. We denote by $B =
\{B(h):h \in \HH\}$ an \textit{isonormal Gaussian process} over
$\HH$, i.e. $B$ is a centered Gaussian family indexed by
the elements of $\HH$ and such that, for every $h_1,h_2\in\HH$,
\begin{equation}\label{isometry}
\E\big[B(h_1)B(h_2)\big]=\langle h_1,h_2\rangle_\HH.
\end{equation}
In what follows, we shall use the notation ${L}^2(B)$ $=$
$L^2(\Omega,\sigma(B),\mathbb P)$. For every $q\geq 1$, we write
$\HH^{\otimes q}$ to indicate the $q$th tensor power of $\HH$; the
symbol $\HH^{\odot q}$ indicates the $q$th \textit{symmetric}
tensor power of $\HH$, equipped with the norm
$\sqrt{q!}\|\cdot\|_{\HH^{\otimes q}}$. We denote by $I_q$ the
isometry between $\HH^{\odot q}$ and the $q$th Wiener chaos of
$X$, which is a linear map satisfying the property
\[
I_q (h^{\otimes q}):= H_q(B(h)), \qquad h^{\otimes q}:=h\otimes \cdots \otimes h\in \HH^{\otimes q} \quad \text{with } \|h\|_{\HH}=1, 
\]
where $H_q$ is the $q$th Hermite polynomial defined in Remark \ref{rem6}. It is well-known (see \cite[Chapter 1]{N}) 
that any random variable $F\in L^2(B)$ admits an orthogonal \textsl{chaotic expansion}:
\begin{equation}\label{ChaosExpansion}
F=\sum_{q=0}^\infty I_q(f_q),
\end{equation}
where $I_0(f_0)=\E[F]$, the series converges in $L^2$ and
the kernels $f_q\in\HH^{\odot q}$, $q\geq1$, are uniquely
determined by $F$. In the particular case where $\HH=L^2(A,\mathcal{A},\mu)$, with $(A,\mathcal{A})$ a measurable
space and $\mu$ a $\sigma$-finite and non-atomic measure, one
has that $\HH^{\odot q}= L^2_s(A^q,\mathcal{A}^{\otimes q},
\mu^{\otimes q})$ is the space of symmetric and square integrable
functions on $A^q$. Moreover, for every $f\in\HH^{\odot q}$,
$I_q(f)$ coincides with the multiple Wiener-It\^o integral (of
order $q$) of $f$ with respect to $B$ (see \cite[Chapter 1]{N}).

Now, we introduce the Malliavin derivative.
Let $\mathcal{S}$ be the set of all smooth cylindrical random
variables of the form $$F = f\big(B(h_1), \ldots,
B(h_n)\big),$$ where $n\geq 1$, $f : \R^n \rightarrow \R$ is a
smooth function with compact support and $h_i\in\HH$. The
Malliavin derivative of $F$  is the element of
$L^2(\Omega, \HH)$ defined as
\[
DF := \sum_{i =1}^n \frac{\partial f}{\partial x_i}\big(B(h_1), \ldots, B(h_n)\big)
h_i.
\]
For instance, $DW(h) = h$ for every $h\in \HH$. We denote by $\mathbb D^{1,2}$  the
closure of $\mathcal{S}$ with respect to the norm $\| \cdot
\|_{1,2}$, defined by the relation
$$\| F\|_{1,2}^2  =  \E[F^2] +
\E[ \| D F\|_{\HH}^2].$$
Note that, if $F$ is equal to a finite sum of multiple Wiener-It\^o
integrals, then $F\in\mathbb D^{1,2}$. The Malliavin
derivative $D$ verifies the following \textit{chain rule}: if
$\varphi:\R^n\rightarrow\R$ is in $C^1_b$ (that is, the
collection of continuously differentiable functions with
bounded partial derivatives) and if $\{F_i\}_{i=1,\ldots,n}$ is a vector of
elements of $\mathbb D^{1,2}$, then
$\varphi(F_1,\ldots,F_n)\in \mathbb D^{1,2}$ and
$$
D\varphi(F_1,\ldots,F_n)=\sum_{i=1}^n
\frac{\partial\varphi}{\partial x_i} (F_1,\ldots, F_n)DF_i.
$$
We denote by $\delta$ the
adjoint of the unbounded operator $D$, also called the \textsl{divergence
operator}. A random element $u \in L^{2}(\Omega, \HH)$
belongs to the domain of $\delta$, noted ${\rm Dom}\delta$, if
and only if it verifies
$$
| \E \langle D F,u\rangle_{\HH}|\leq c_u\,\|F\|_{L^2}\quad\mbox{for any }F\in{\mathcal{S}},
$$
where $c_u$ is a constant depending only on $u$. If $u \in
{\rm Dom}\delta $, then the random variable $ \delta(u)$ is
defined by the duality relationship (sometimes called `integration by
parts formula'):
\begin{equation}\label{ipp}
\E [F \delta(u)]=  \E \langle D F, u \rangle_{\HH},
\end{equation}
which holds for every $F \in \mathbb D^{1,2}$. An immediate consequence of \eqref{ipp}
is the following identity
\begin{align} \label{ipp2}
\delta (Fu) = F\delta (u) - \langle DF,u \rangle_{\HH},
\end{align}
which holds for all $F \in \mathbb D^{1,2}$ and $u \in {\rm Dom}\delta $ such that $Fu\in {\rm Dom}\delta$.

\begin{rem} \label{hilbert} \rm
When $\HH=L^2([0,T],dx)$, which is the most basic example,  
then the isonormal Gaussian process $B$ is a standard Brownian motion on $[0,T]$.
In this case we have 
\[
\delta(h) = \int_0^T h_s dB_s, \qquad \forall h\in L^2([0,T],dx). 
\]
The divergence operator $\delta$ is often called Skorohod integral. One can show that for a stochastic process 
$u \in {\rm Dom}\delta$, the Skorohod integral $\delta (u)$ and the It\^o integral $\int_0^T u_s dB_s$ coincide
whenever the latter is well defined.

In case $\HH=L^2([0,T],dx)$ the Malliavin derivative $D$ and the divergence operator $\delta$ can be computed directly using chaos expansion.
Indeed,  
the derivative of a random variable
$F$ as in (\ref{ChaosExpansion}) can be identified with the
element of $L^2([0,T]\times \Omega)$ given by
\begin{equation}\label{dtf}
D_a F=\sum_{q=1}^\infty qI_{q-1}\big(f_q(\cdot,a)\big), \quad a
\in [0,T].
\end{equation}
On the other hand, for any $u \in L^2([0,T]\times \Omega)$ there exists a chaos decomposition
\[
u_s = \sum_{q=0}^{\infty} I_q (f_q(\cdot, s)), \qquad  f_q(\cdot, s) \in  L^2_s([0,T]^q\times \Omega).
\] 
Let $\tilde{f}_q\in L^2_s([0,T]^{q+1}\times \Omega)$ denote the symmetrization of $f_q(\cdot, \cdot)$. Then 
the element $\delta(u)$ can be written in terms of chaotic decomposition as
\[
\delta(u) = \sum_{q=0}^{\infty} I_{q+1} (\tilde{f}_q).
\] 
\qed
\end{rem}
Now, we start introducing the definition of the integral $\int_0^t Z_s dX_s$, where the ambit process $X$ is driven by a Brownian motion $W$.
The following exposition is related to a seminal work \cite{AMN}, where the integration with respect to Gaussian processes has been investigated. 
Throughout this subsection we assume that
\begin{align} \label{l2cond}
\E\left[ \int_0^t g^2(t,s)\sigma_s^2 ds \right]<\infty 
\end{align} 
and $\mathcal F=\sigma (W_t:~t\in [0,T])$. The following definition is due to \cite[Section 4]{BBPV}.

\begin{defi} \label{def2} \rm
Assume that, for any $s\geq 0$, the function $t\mapsto g(t,s)$ has bounded variation on the interval $[t,v]$ for all $0\leq s< t<v<\infty $.
We say that the process $Z$ belongs to the class $I^X(0,T)$ when the following conditions are satisfied:
\begin{itemize}
\item[(i)] For any $s\in [0,T]$ the process $(Z_u - Z_s)_{u\in (s,T]}$ is integrable with respect to $g(du,s)$.
\item[(ii)] Define the operator 
\[
\mathcal{K}_g(h)(t,s):= h(s) g(t,s)+ \int_s^t (h(u)-h(s)) g(du,s).
\]
The process $s\mapsto \mathcal{K}_g(Z)(T,s)\sigma_s 1_{[0,T]}(s)$ belongs to ${\rm Dom}\delta$.
\item[(iii)] $\mathcal{K}_g(Z)(T,s)$ is Malliavin differentiable with respect to $D_s$ with $s\in [0,T]$, such that
the mapping $s\mapsto D_s[\mathcal{K}_g(Z)(T,s)]\sigma_s$ is Lebesgue integrable on $[0,T]$. 
\end{itemize}
When $Z\in I^X(0,T)$ we define
\[
\int_0^T Z_s dX_s:= \delta \left( \mathcal{K}_g(Z)(T,s)\sigma_s 1_{[0,T]}(s) \right) + \int_0^T D_s[\mathcal{K}_g(Z)(T,s)]\sigma_s ds.
\]
\end{defi}
We remark that the proposed definition is linear in the integrand. It also holds that
\[
\int_0^T YZ_s dX_s = Y \int_0^T Z_s dX_s
\]
for any bounded random variable $Y$ such that $Z, YZ\in I^X(0,T)$. We refer to \cite{BBPV} for further properties and applications. 

The operator $\mathcal{K}$ has been introduced in \cite{AMN}. The intuition behind Definition \ref{def2} is explained by the following
heuristic derivation. Using classical integration by parts formula and \eqref{ipp2} we conclude that
\begin{align*}
\int_0^T Z_s dX_s &= Z_T X_T - \int_0^T \frac{dZ_u}{du} \left(\int_0^u g(u,s) \sigma_s dW_s\right) du \\[1.5 ex]
&= Z_T X_T - \int_0^T \delta\left( \frac{dZ_u}{du}  g(u,s) \sigma_s 1_{[0,T]}(s) \right) du \\[1.5 ex]
&- \int_0^T \int_0^u D_s \left[ \frac{dZ_u}{du} \right] g(u,s) \sigma_s ds du.
\end{align*}  
Next, the stochastic Fubini theorem applied to the last two quantities implies the identity
\begin{align*}
\int_0^T Z_s dX_s &= Z_T X_T - \delta \left( \sigma_s \int_s^T g(u,s) \frac{dZ_u}{du} du 1_{[0,T]} (s)\right) \\[1.5 ex]
&- \int_0^T  D_s \left[ \int_s^T g(u,s) \frac{dZ_u}{du} du\right] \sigma_s ds.
\end{align*} 
Similarly, we deduce that 
\[
Z_T X_T = \delta \left( Z_T g(T,s)\sigma_s  1_{[0,T]} (s)\right) + \int_0^T  D_s \left[ Z_T\right] g(T,s) \sigma_s ds.
\]
Thus, putting things together and applying the classical integration by parts formula once again we obtain the heuristic formula
\begin{align*}
\int_0^T Z_s dX_s &= \delta \left( \sigma_s \left[Z_Tg(T,s) - \int_s^T g(u,s) \frac{dZ_u}{du} du  \right] 1_{[0,T]} (s) \right) \\[1.5 ex]
&- \int_0^T  D_s \left[Z_T g(T,s) -\int_s^T g(u,s) \frac{dZ_u}{du} du\right] \sigma_s ds \\[1.5 ex]
& = \delta \left( \mathcal{K}_g(Z)(T,s)\sigma_s 1_{[0,T]}(s) \right) + \int_0^T D_s[\mathcal{K}_g(Z)(T,s)]\sigma_s ds.
\end{align*}
This explains the intuition behind Definition \ref{def2}.

\begin{rem} \rm
In a recent work \cite{BS13} the integration concept has been extended to the class of Hilbert-valued processes. \qed
\end{rem}

\subsection{Integration with respect to ambit processes driven by pure jump L\'evy motion}
In this subsection we will introduce the definition of the integral 
\begin{align*} 
\int_0^t Z_s dX_s,
\end{align*} 
where the ambit process $X$ is driven by a square integrable pure jump L\'evy motion $L$ with characteristic triplet $(0,0,\nu)$. 
We assume that the condition \eqref{l2cond} holds and $\mathcal F =\sigma(L_t:~t\in [0,T])$. 

The definition of the stochastic integral proposed in \cite{BBPV} relies again on Malliavin calculus. However, in contrast to the Gaussian space,
there exist different variations of Malliavin calculus for Poisson random measures. Here we follow an approach described in  \cite{NOP}.
We deal with the Hilbert space $\HH=L^2([0,T],dx)$. Let $N(dt,dz)$ denote the Poisson random measure on $[0,T] \times (\R\setminus \{0\})$
associated with $L$ and $\widetilde N(dt,dz) = N(dt,dz) - dt\nu (dz)$ the compensated Poisson random measure. We have
\begin{align*}
L_t = \int_0^t \int_{\R\setminus \{0\}} z \widetilde N(dt,dz). 
\end{align*} 
As in the Gaussian case there exists an orthogonal chaos decomposition of the type \eqref{ChaosExpansion} in terms of multiple integrals. 
For any $f\in L^2_s (([0,T] \times (\R\setminus \{0\}))^q)$, we introduce the $q$th order multiple integral of $f$ with respect to $\widetilde N(dt,dz)$
via
\begin{align*}
I_q(f) = q! \int_0^T \int_{\R\setminus \{0\}} \cdots \int_0^{t_2-} \int_{\R\setminus \{0\}}  f(t_1,z_1, \ldots, t_q,z_q) \widetilde N(dt_1,dz_1)\ldots
N(dt_q,dz_q). 
\end{align*} 
Then, for any random variable $F\in L^2(\Omega , \mathcal F, \mathbb P)$, there exists a unique sequence of symmetric functions
$(f_q)_{q\geq 0}$ with $f_q\in L^2_s (([0,T] \times (\R\setminus \{0\}))^q)$ such that
\begin{align} \label{chaos2}
F= \sum_{q=0}^{\infty } I_q(f_q),
\end{align} 
which is obviously an analogue of \eqref{ChaosExpansion}. Furthermore, it holds that
\[
\E[F^2] = \sum_{q=0}^{\infty } q! \|f_q\|_{q,\nu}^2,
\]
where the norm $\|f_q\|_{q,\nu}^2$ is defined by
\[
\|f_q\|_{q,\nu}^2 := \int_{([0,T] \times (\R\setminus \{0\}))^q} f^2(t_1,z_1, \ldots, t_q,z_q) dt_1\nu (dz_1) \ldots dt_q\nu (dz_q). 
\]
Similarly to the exposition of Remark \ref{hilbert}, the Malliavin derivative $D$ and the divergence operator $\delta$
are introduced using the above chaos representation. We say that a random variable $F$  with chaos decomposition \eqref{chaos2}
belongs to the space $\mathbb D^{1,2}$
whenever the condition 
\[
\sum_{q=0}^{\infty } q q! \|f_q\|_{q,\nu}^2<\infty 
\]   
holds. Whenever $F\in \mathbb D^{1,2}$ we define
\[
D_{t,z} F := \sum_{q=1}^{\infty } q I_{q-1}(f_q(\cdot,t,z)).
\] 
Now, we say that that a random field $u\in L^2([0,T] \times \R\setminus \{0\} \times \Omega )$ belongs to the domain 
of the divergence operator $\delta$ (${\rm Dom}\delta$) when
\[
\left| \E \left[ \int_0^T \int_{\R\setminus \{0\}} u(t,z) D_{t,z} F \nu (dz)dt  \right] \right| \leq c_u \|F\|_{L^2}
\]
for all $F\in \mathbb D^{1,2}$. Whenever $u\in {\rm Dom}\delta$ the element $\delta (u)$ is uniquely characterized via the identity
\[
\E[F\delta (u)]= \E \left[ \int_0^T \int_{\R\setminus \{0\}} u(t,z) D_{t,z} F \nu (dz)dt  \right] \qquad \forall F\in \mathbb D^{1,2}, 
\]
which is an integration by parts formula (cf. \eqref{ipp}). An immediate consequence of the integration by parts formula
is the following equation:
\begin{align}
F \delta (u) = \delta (u(F+DF)) + \int_0^T \int_{\R\setminus \{0\}} u(t,z) D_{t,z} F \nu (dz)dt,
\end{align}
which holds for any $F\in \mathbb D^{1,2}$, $u\in {\rm Dom}\delta$ such that $u(F+DF)\in {\rm Dom}\delta$. Notice the appearance
of the term $uDF$ on the right hand side that is absent in the Gaussian case (cf. \eqref{ipp2}).

Now, we proceed with the introduction of the stochastic integral.
Its definition in the pure jump L\'evy framework is essentially analogous to the Gaussian case. We refer again to \cite[Section 4]{BBPV}
for a more detailed exposition and the intuition behind this definition.    

\begin{defi} \label{def3} \rm
Assume that, for any $s\geq 0$, the function $t\mapsto g(t,s)$ has bounded variation on the interval $[t,T]$.
We say that the process $Z$ belongs to the class $I^X(0,T)$ when the following conditions are satisfied:
\begin{itemize}
\item[(i)] For any $s\in [0,T]$ the process $(Z_u - Z_s)_{u\in (s,T]}$ is integrable with respect to $g(du,s)$.
\item[(ii)] 
The process $(s,z) \mapsto z(\mathcal{K}_g(Z)(T,s)+ D_{s,z}[\mathcal{K}_g(Z)(T,s)])\sigma_s 1_{[0,T]}(s)$ belongs to ${\rm Dom}\delta$.
\item[(iii)] $\mathcal{K}_g(Z)(T,s)$ is Malliavin differentiable with respect to $D_{s,z}$ with $(s,z)\in [0,T] \times \R$, such that
the mapping $(s,z)\mapsto z D_{s,z}[\mathcal{K}_g(Z)(T,s)]\sigma_s$ is $\nu(dz)dt$-integrable. 
\end{itemize}
When $Z\in I^X(0,T)$ we define
\begin{align*} 
\int_0^T Z_s dX_s&:= \delta \left(z(\mathcal{K}_g(Z)(T,s)+ D_{s,z}[\mathcal{K}_g(Z)(T,s)])\sigma_s 1_{[0,T]}(s) \right) \\[1.5 ex]
&+ \int_0^T \int_{\R} z D_{s,z}[\mathcal{K}_g(Z)(T,s)]\sigma_s \nu(dz)ds.
\end{align*}
\end{defi}

\section{Limit theory for high frequency observations of ambit fields} \label{sec5}
\setcounter{equation}{0}
\renewcommand{\theequation}{\thesection.\arabic{equation}}
In this section we will review the asymptotic results for power variation of L\'evy semi-stationary processes (\textit{LSS}) 
without drift, i.e.
\[
Y_{t}=\mu +\int_{-\infty}^t g(t-s)\sigma_{s} L(\mathrm{d}s). 
\]
We will see that the limit theory heavily depends on whether the driving L\'evy motion is a Brownian motion or a pure jump
process. Furthermore, the structure of the weight function $g$ plays an important role as we will see below. More precisely,
the \textit{singularity points} of $g$ determine the type of the limit.

In what follows we assume that the underlying observations of L\'evy semi-stationary process $Y$ are
\[
Y_0, Y_{\Delta_n}, Y_{2\Delta_n}, \ldots, Y_{\Delta_n [t/\Delta_n]}
\] 
with $\Delta_n\rightarrow 0$ and $t$ fixed. In other words, we are in the infill asymptotics setting. For statistical purposes
we introduce $k$th order differences $\Delta_{i,k}^{n} Y$ of $Y$ defined via 
\begin{align} \label{filter}
\Delta_{i,k}^{n} Y:= \sum_{j=0}^k (-1)^j \binom{k}{j} Y_{(i-j)\Delta_n}.
\end{align}
For instance, 
\[
\Delta_{i,1}^n Y = Y_{i\Delta_n}-Y_{(i-1)\Delta_n} \quad \textrm{and} \quad \Delta_{i,2}^n Y = Y_{i\Delta_n}-2Y_{(i-1)\Delta_n}+Y_{(i-2)\Delta_n}.
\] 
The power variation of $k$th order differences of $Y$ is given by the statistic
\begin{align} \label{powervar}
V(Y,p,k;\Delta_n)_t &:= \sum_{i=k}^{[t/\Delta_n]} |\Delta_{i,k}^{n} Y|^p. 
\end{align} 
In the following we will study the asymptotic behaviour of the functional $V(Y,p,k;\Delta_n)_t$.

\subsection{\textit{LSS} processes driven by Brownian motion}
In this section we consider the case of \textit{Brownian semi-stationary processes} given via 
\[
Y_{t}=\mu +\int_{-\infty}^t g(t-s)\sigma_{s} W(\mathrm{d}s),
\]
defined on a filtered probability space $(\Omega ,\mathcal{F},(\mathcal{F}_{t})_{t\in \mathbb{R}},\mathbb{P})$.
The following asymptotic results have been investigated in a series of papers \cite{BCP11, BCP13, CHPP13, GP}.
We also refer to a related work \cite{BCP09,CNW}, where the power variation of integral processes as defined 
in Remark \ref{rem8} has been studied.  

In the model introduced above $W$ is an $(\mathcal{F}_{t})_{t\in \mathbb{R}}$-adapted white noise on $\mathbb{R}$, 
$g:\R \rightarrow \mathbb{R}$ is a deterministic weight function satisfying
$g(t)=0$ for $t\leq 0$ and $g\in \mathbb{L}^{2}(\mathbb{R})$. The intermittency  process $\sigma $ 
is assumed to be an \ $(\mathcal{F}_{t})_{t\in \mathbb{R}}$-adapted c\`{a}dl\`{a}g process. 
The finiteness of the process $X$ is guaranteed by the condition
\begin{align} \label{fincon}
\int_{-\infty }^{t}g^2(t-s)\sigma^2_{s} ds<\infty \quad \text{almost surely}, 
\end{align}
for any $t\in \R$, which we assume from now on. As pointed out in \cite{BCP11,BCP13} and also 
briefly discussed in Remark \ref{rem8}, the \textit{Gaussian core} $G$ is crucial for understanding the fine structure of $Y$. 
The process $G=(G_t)_{t\in \R}$ is a zero-mean stationary Gaussian process given by
\begin{align} \label{g}
G_t: = \int_{-\infty}^t g(t-s)  W(ds), \qquad t\in \R.
\end{align}  
We remark that $G_t<\infty$ since $g\in \mathbb{L}^2(\R)$. A straightforward computation
shows that the correlation kernel $r$ of $G$ has the form
\[
r(t) = \frac{\int_0^\infty g(u)g(u+t) du}{\|g \|_{\mathbb L^2(\R)}^2}, \qquad t\geq 0.
\]
Another important quantity for the asymptotic theory is the variogram  $R$, i.e.
\begin{align} \label{R}
R(t):= \E[(G_{t+s} - G_s)^2] = 2 \|g \|_{\mathbb L^2(\R)}^2 (1-r(t)), \qquad 
\tau_{k}(\Delta_n):=\sqrt{\E[(\Delta_{i,k}^{n} G)^2]}.
\end{align} 
The quantity $\tau_{k}(\Delta_n)$ will appear as a proper scaling in the law of large numbers
for the statistic $V(Y,p,k;\Delta_n)$ introduced at \eqref{powervar}.

As mentioned above the set of singularity points $0=\theta_0<\theta_1<\cdots<\theta_l<\infty $ of $g$ will determine
the limit theory for the power variation $V(Y,p,k;\Delta_n)$. Let 
$\alpha_0,\ldots ,\alpha_l\in (-1/2,0)\cup (0,1/2)$ be given real numbers.
For any function $h\in C^m (\R)$, $h^{(m)}$ denotes the $m$-th derivative of $h$. Recall that $k\geq 1$ stands for the order of the
filter defined in \eqref{filter}. We introduce the following set of assumptions. \\ \\
(A): For $\delta<\frac 12 \min_{1\leq i\leq l}(\theta_i-\theta_{i-1})$ it holds that \\ \\
(i) $g(x)= x^{\alpha_0}  f_0(x)$ for $x\in (0, \delta)$ and $g(x)= |x-\theta_l|^{\alpha_l}  f_l(x)$
for $x\in (\theta_l-\delta, \theta_l)\cup (\theta_l, \infty )$. \\ \\
(ii) $g(x)= |x-\theta_i|^{\alpha_i}  f_i(x)$ for $x\in (\theta_i-\delta, \theta_i)\cup (\theta_i, \theta_i +\delta)$,
$i=1,\ldots,l-1$. \\ \\
(iii) $g(\theta_i)=0$, 
$f_i\in C^k\left( (\theta_i-\delta, \theta_i +\delta)\right)$ and $f_i(\theta_i)\not =0$ for $i=0,\ldots, l$.\\ \\
(iv) $g\in C^k(\R \setminus \{\theta_0, \ldots, \theta_l\})$ and $g^{(k)}\in \mathbb L^2 \left(
\R \setminus \cup_{i=0}^l (\theta_i-\delta, \theta_i+\delta) \right)$. \\ \\
(v) For any $t>0$
\begin{align} \label{F}
F_t= \int_{\theta_l+1}^\infty g^{(k)}(s)^2 \sigma_{t-s}^2 ds <\infty.  
\end{align} 
We also set
\begin{align} \label{min}
\alpha := \min\{\alpha_0,\ldots ,\alpha_l \}, \qquad \mathcal{A}:=\{0\leq i\leq l:~\alpha_i=\alpha\}.
\end{align}
The points $\theta_0,\ldots,\theta_l$ are singularities of $g$ in the sense that $g^{(k)}$ is not square integrable 
around these points, because $\alpha_0,\ldots ,\alpha_l\in (-1/2,0)\cup (0,1/2)$ and conditions (A)(i)-(iii) hold. 
Condition (A)(iv) indicates that $g$ exhibits no further singularities.

\begin{rem} \label{rem9} \rm
According to the discussion of Section \ref{sec2.3}, the Brownian semi-stationary processes $Y$ (or even
the Gaussian core $G$) is not a semimartingale, since $g'\not \in L^2(\R_{\geq 0})$
due to the presence of the singularity points $\theta_0,\ldots,\theta_l$. For this reason we can not rely
on limit theory for power variations of continuous semimartingales investigated in e.g. \cite{BGJPS06,J08}.
Although some of the asymptotic results look similar to the semimartingale case, the methodology behind
the proof is completely different. The main steps of the proof are based on methods of Malliavin calculus
developed in e.g. \cite{NP05,PT05} and on Bernstein's blocking technique. \qed 
\end{rem}  
The limit theory for the power variation is quite different according to whether we have a single singularity at, say, $\theta_0=0$
(i.e. $l=0$) or multiple singularity points. Hence, we will treat the corresponding results separately.

\subsubsection{The case $l=0$}
The theory presented in this section is mainly investigated in \cite{BCP11,BCP13}. Below we will intensively use
the concept of stable convergence, which is originally due to R\'enyi \cite{R63}. 
We say that a sequence of processes $X^n$ converges stably in law to
a  process $X$, where $X$ is defined on an extension $(\Omega', \mathcal{F}', \mathbb P')$ 
of the original probability $(\Omega, \mathcal{F}, \mathbb P)$, in the space $\mathbb D([0,T])$
equipped with the uniform topology ($X^n \stab X$) if and only if
\begin{equation*} 
\lim_{n\rightarrow \infty} \E[f(X^n) Z] = \E'[f(X)Z] 
\end{equation*} 
for any bounded and continuous function $f: \mathbb D([0,T]) \rightarrow \R$ and any bounded $\mathcal F$-measurable
random variable $Z$.  We refer to \cite{AE78}, \cite{JS02} or \cite{R63} for a detailed study of stable convergence. 
Note that stable convergence is a stronger mode of convergence than weak convergence, but it is weaker that 
u.c.p. convergence. 

The following theorem has been shown in \cite[Theorems 1 and 2]{BCP13}.

\begin{theo} \label{th5}
Assume that condition (A) holds.
\begin{itemize}
\item[(i)] We obtain that 
\begin{align} \label{lln}
\Delta_n \tau_{k}(\Delta_n)^{-p}V(Y,p,k;\Delta_n)_t \ucp V(Y,p)_t:= m_p \int_0^t |\sigma_s|^p ds,
\end{align}
where the power variation $V(Y,p,k;\Delta_n)_t$ is defined at \eqref{powervar} and the constant $m_p$ is given
by \eqref{mp}.
\item[(ii)] Assume that the intermittency process $\sigma$ is H\"older continuous of order $\gamma \in (0,1)$
and $\gamma (p\wedge 1)>1/2$. When $k=1$ we further assume that $\alpha \in (-1/2,0)$. Then we obtain the stable convergence
\begin{align} \label{clt}
\Delta_n^{-1/2} \Big(\Delta_n \tau_{k}(\Delta_n)^{-p}V(Y,p,k;\Delta_n)_t - V(Y,p)_t \Big) 
\stab  \lambda \int_0^t |\sigma_s|^p  ~dB_s
\end{align}
on $\mathbb D([0,T])$ equipped with the uniform topology, where 
$B$ is a  Brownian motion that is defined on an
extension of the original probability space $(\Omega, \mathcal{F}, \mathbb P)$ and is independent of $\mathcal{F}$, 
and the constant $\lambda$ is  given by
\begin{align} 
\lambda^2= \lim_{n\rightarrow \infty } \De_n^{-1} \mathrm{var} \Big( \Delta_n^{1-pH} V(B^H,p,k;\Delta_n)_1\Big), \qquad 
\end{align}
with $B^H$ being a fractional Brownian motion with Hurst parameter $H=\alpha + 1/2$.
\end{itemize}
\end{theo}

\begin{rem} \label{rem10} \rm
The appearance of the fractional Brownian motion in the definition of the constant $\lambda^2$ is not surprising 
given the discussion of fine properties of the Gaussian core $G$ in Section \ref{sec2.3}. In case $k=1$ the factor $\lambda^2$
coincides with the quantity $v_p$ defined in Remark \ref{rem6}. We also remark that the validity region of the central limit theorem
in \eqref{clt} in the case $k=1$ ($\alpha \in (-1/2,0)$) is smaller than the region $H=\alpha +1/2\in (0,3/4)$ described in 
Remark \ref{rem6}. This is due to a bias problem, which appears in the context of Brownian semi-stationary processes. \qed     
\end{rem}
Notice that the asymptotic result at \eqref{lln} and \eqref{clt} are not feasible from the statistical point of view, since
the scaling $\tau_{k}(\Delta_n)$ depends on the unknown weight function $g$. Nevertheless, Theorem \ref{th5} is useful for statistical
applications. Our first example is the estimation of the \textit{smoothness parameter} $\alpha$. Under mild conditions
on the intermittency process $\sigma$, the Brownian semi-stationary process $Y$ (as well as its Gaussian core $G$) has 
H\"older continuous paths of any order smaller than $H=\alpha +1/2\in (0,1)$. In turbulence the smoothness parameter $\alpha$
is related to the so called \textit{Kolmogorov's $2/3$-law}, which predicts that
\[
\E[(X_{t+\Delta} -X_t)^2]\propto \Delta^{2/3}, 
\]  
or in other words $\alpha \approx -1/6$, which holds for a certain range of frequencies $ \Delta $. 
Hence, estimation of the parameter $\alpha$ is extremely important. 

A typical model for the weight function $g$ is the Gamma kernel given via
\[
g(x)=x^{\alpha} \exp(-cx), \qquad c>0, \alpha \in (-1/2,0)\cup (0,1/2),  
\] 
which obviously satisfies the assumption (A)(i)-(iv) with $l=0$. An application of the law of large numbers at 
\eqref{lln} for a fixed $t>0$ gives 
\begin{align*}
S_n:=\frac{V(Y,p,k;2\Delta_n)_t}{V(Y,p,k;\Delta_n)_t} \toop 2^{\frac{(2\alpha +1)p}{2}},
\end{align*}      
since $\tau_{k}(2\Delta_n)^2/ \tau_{k}(\Delta_n)^2\rightarrow 2^{2\alpha +1}$. The latter is due
to $\tau_{k}(\Delta_n)^2 \sim \Delta_n^{2\alpha +1}$, which follows from the fact that the Gaussian core $G$
and the fractional Brownian motion $B^H$ with Hurst parameter $H=\alpha +1/2$ have the same
small scale behaviour. Thus,  a consistent estimator of $\alpha$ is given via
\begin{align} \label{hatalpha}
\widehat{\alpha}_n=\frac 12 \left( \frac{2\log_2 S_n}{p}  -1
\right) \toop \alpha ,
\end{align}
where $\log_2$ denotes the logarithm at basis $2$. Note that the estimator $\widehat{\alpha}_n$ is feasible, i.e. it does not
depend on the unknown scaling $\tau_{k}(\Delta_n)$. One may also deduce a standard feasible central limit theorem for
$\widehat{\alpha}_n$ as it was shown in \cite{BCP11,BCP13,CHPP13}, and thus obtain asymptotic confidence regions for
the smoothness parameter $\alpha$. We also refer to \cite{CHPP13} for empirical implementation of this estimation method 
to turbulence data. 

Another useful application of Theorem \ref{th5} is the estimation of the \textit{relative intermittency}, which is defined as
\[
RI_t:= \frac{\int_0^t \sigma_s^2 ds}{\int_0^T \sigma_s^2 ds}, \qquad t\leq T,
\] 
where $T>0$ is a fixed time. While the intermittency process $\sigma$ is not identifiable when no structural assumption 
on $g$ are imposed, the relative intermittency $RI_t$ is easy to estimate. Indeed, the convergence in \eqref{lln} 
immediately implies that
\[
\widehat{RI}_t^n:=\frac{V(Y,p,k;\Delta_n)_t}{V(Y,p,k;\Delta_n)_T} \toop RI_t. 
\] 
We refer to \cite{BPS} for the limit theory and physical applications of the statistic $\widehat{RI}_t^n$.

\subsubsection{The case $l\geq 1$} 
The limit theory for the case $l\geq 1$ appears to be more complex. The asymptotic results presented below
have been investigated in \cite{GP}.

First of all, we need to introduce some notations. Recall that $k\in \mathbb N$ denotes the order of increments. 
The $k$-th order filter associated with $g$ is introduced via
\begin{align} \label{filterg}
\Delta_{k}^{n} g(x):= \sum_{j=0}^k (-1)^j \binom{k}{j} g(x-j\Delta_n), \qquad x\in \R.
\end{align}
There is a straightforward relationship between the scaling quantity $\tau_{k}(\Delta_n)$
defined at \eqref{R} and the function $\Delta_{k}^{n,} g$, namely
\[
\tau_{k}(\Delta_n)^2 = \|\Delta_{k}^{n} g\|_{\mathbb{L}^2(\R)}^2.
\] 
Now, we define the concentration measure associated with $\Delta_{k}^{n} g$:
\begin{align} \label{pink}
\pi_{n,k}(A) := \frac{\int_A (\Delta_{k}^{n} g(x))^2 dx}{\|\Delta_{k}^{n} g\|_{\mathbb{L}^2(\R)}^2},
\qquad A\in \mathcal{B}(\R).
\end{align}
Observe that $\pi_{n,k}$ is a probability measure. Its asymptotic behaviour determines the law of large numbers
for the power variation. In order to identify the limit of $\pi_{n,k}$, we define the following functions
\begin{align} \label{hdef}
h_0 (x)&:= f_0(\theta_0) \sum_{j=0}^k (-1)^j \binom{k}{j} (x-j)_{+}^{\alpha_0}, \\[1.5 ex]
h_i (x)&:= f_i(\theta_i) \sum_{j=0}^k (-1)^j \binom{k}{j} |x-j|^{\alpha_i}, \qquad i=1,\ldots, l, \nonumber
\end{align}
where $x_+:=\max\{x,0\}$. The following results determines the asymptotic behaviour of the power variation 
$V(Y,p,k;\Delta_n)_t$ for $p=2$ and $l\geq 1$. We refer to \cite[Proposition 3.1, Theorems 3.2 and 3.3]{GP} for further details.

\begin{theo} \label{th6}
Assume that the condition (A) holds.
\begin{itemize}
\item[(i)] It holds that 
\begin{align*}
\pi_{n,k} \schw \pi_k,
\end{align*}
for any $k\geq 1$, where the probability measure $\pi_k$ is given as
\begin{align} \label{pik}
\text{supp}(\pi_k)= \{\theta_i\}_{i\in \mathcal{A}}, \qquad 
\pi_k (\theta_i)= \frac{\|h_i\|_{\mathbb{L}^2(\R)}^2 1_{i\in \mathcal{A}}}{\sum_{j=0}^l \|h_j\|_{\mathbb{L}^2(\R)}^2
1_{j\in \mathcal{A}}},
\end{align}
where the set $\mathcal{A}$ has been defined at \eqref{min}.
\item[(ii)] We obtain the convergence
\begin{align} \label{LLN}
\frac{\Delta_n}{\tau_{k}(\Delta_n)^2} V(Y,2,k;\Delta_n)_t \ucp   QV(Y,k)_t:= 
\int_0^\infty \left( \int_{-\theta}^{t-\theta } \sigma_s^2 ds\right) \pi_k (d\theta).
\end{align}
\item[(iii)] Assume that the intermittency process $\sigma$ is H\"older continuous of order
$\gamma >1/2$. When $k=1$ we further assume that $\alpha_j \in (-\frac 12,0)$ for all $0\leq j\leq l$. Then, under condition
\begin{align} \label{robustness}
\alpha_i -\alpha > 1/4 \qquad \text{for all } i \not \in \mathcal{A},  
\end{align} 
we obtain the stable convergence
\begin{align} \label{CLT} 
\Delta_n^{-1/2} \left(\frac{\Delta_n}{\tau_{k}(\Delta_n)^2} V(Y,2,k;\Delta_n)_t -   QV(X,k)_t \right) 
\stab   \int_0^t  v_s^{1/2} dB_s
\end{align}
on $\mathbb D([0, \min_{1\leq j\leq l}(\theta_j-\theta_{j-1})])$ equipped with the uniform topology, where 
$B$ is a  Brownian motion, independent of $\mathcal F$, defined on an
extension of the original probability space $(\Omega, \mathcal{F}, \mathbb P)$.
The stochastic process $v$ is given by 
\begin{align} \label{asyvar}
v_s= \Lambda_k \left(\int_0^{\infty} \sigma^2_{s-\theta} \pi_k (d\theta)\right)^2 ,
\end{align}
where $\Lambda_k$ is  defined by
\[
\Lambda_k = \lim_{n\rightarrow \infty} \De_n^{-1} \mathrm{var} \Big( 
\frac{\Delta_n}{\hat{\tau}_{k}(\Delta_n)^2} V(B^H,2,k;\Delta_n)_1\Big)
\]
with $B^H$ being a fractional Brownian motion with Hurst parameter $H=\alpha + 1/2$ and 
$\hat{\tau}_{k}(\Delta_n)^2:=\E[(\Delta_{i,k}^{n} B^H)^2]$.
\end{itemize}
\end{theo} 
In order to explain the mathematical intuition behind the results of Theorem \ref{th6} we present some remarks.

\begin{rem} \label{rem11} \rm
We notice that $\text{supp}(\pi_k)= \{\theta_i\}_{i\in \mathcal{A}}$, which means that only those singularity
points contribute to the limit, which correspond to the minimal index $\alpha$. This fact is not surprising 
from the statistical point of view, since a process with the roughest path always dominates when 
considering a power variation. 

When $l=0$ it holds that $\pi_{n,k} \schw \delta_{\{0\}}$, hence the convergence in \eqref{lln} is a particular
case of \eqref{LLN}. Otherwise, the limiting measure $\pi_k$ is a discrete probability measure. It is an open
problem whether the result of \eqref{LLN} can be deduced for a continuous probability measure $\pi_k$. \qed    
\end{rem}

\begin{rem} \label{rem12} \rm
Although the singularity points $\theta_i$ with $i \not \in \mathcal{A}$ do not contribute to the limit at 
\eqref{LLN}, they cause a certain bias, which might explode in the central limit theorem. Condition  
\eqref{robustness} ensures that it does not happen. We also remark that the functional stable convergence at 
\eqref{CLT} does not hold on any interval $[0,T]$, but just on $[0, \min_{1\leq j\leq l}(\theta_j-\theta_{j-1})]$.
One may still show a stable central limit theorem with an $\mathcal F$-conditional Gaussian
process as the limit on a larger interval, 
but only when $\theta_j-\theta_{j-1}\in \mathbb N$ for all $j$, since otherwise the covariance structure
of the original statistic does not converge. \qed
\end{rem}
Notice that the minimal parameter $\alpha$ defined at \eqref{min} still determines the H\"older continuity of the process
$Y$ (and the Gaussian core $G$). The estimator $\widehat{\alpha}_n$ defined at \eqref{hatalpha}  remains consistent, i.e.
\[
\widehat{\alpha}_n=\frac 12 \left( \log_2 \frac{V(Y,2,k;2\Delta_n)_t}{V(Y,2,k;\Delta_n)_t}  -1
\right) \toop \alpha.
\]  
One may also construct a standardized version of the statistic $\widehat{\alpha}_n$, which satisfies a standard central limit theorem
(see \cite[Section 4]{GP} for a detailed exposition). But in this case the time $t<\min_{1\leq j\leq l}(\theta_j-\theta_{j-1})$
must be used, which requires the knowledge of singularity points $\theta_i$.

\begin{rem} \label{rem13} \rm
For potential applications in turbulence the asymptotic results need to be extended to ambit fields $X$ driven by a Gaussian random measure, i.e.
\[
X_{t}(x)=\mu +\int_{A_{t}(x)}g(t,s,x,\xi )\sigma_{s}(\xi )W(\mathrm{d}s,\mathrm{d}\xi ), \qquad t\geq 0, x\in \R^3,
\] 
where $W$ is a Gaussian random measure. This type of high frequency limit theory has not been yet investigated neither
for observations of $X$ on a grid in time-space nor for observations of $X$ along a curve. In the multiparameter setting
there exist a related work on generalized variation of fractional Brownian sheet (see e.g. \cite{PR}) and integral 
processes (see e.g. \cite{Pa,Re}).  
\qed 
\end{rem}

\subsection{\textit{LSS} processes driven by a pure jump L\'evy motion}
In this section we will mainly study the power variation of a L\'evy moving average process defined via
\[
Y_{t}=\mu +\int_{-\infty}^t g(t-s) dL_s,
\] 
where $L$ is a two sided pure jump L\'evy motion without drift. Notice that this is a subclass of \textit{LSS} processes
with $\sigma =1$, and hence it plays a similar role as the Gaussian core $G$ defined at \eqref{g}. The asymptotic theory for power
variation of $Y$ is likely to transfer to limit theory for general \textit{LSS} processes as it was indicated in Remark \ref{rem8}.
In contrast to Gaussian moving averages, little is known about power variation of  L\'evy moving average processes driven by a pure
jump L\'evy motion. Below we present some recent results from \cite{BLP}, which completely determine the first order structure
of power variation of $Y$. We need to impose a somewhat similar set of assumptions as presented in (A) for case $l=0$ (in particular, they insure
the existence of $Y_t$, cf. Remark \ref{rem4}). \\ \\
(A'): It holds that \\ \\
(i) $g(x)= x^{\alpha}  f(x)$ for $x\geq 0$  and $g(x)=0$ for $x<0$ with $\alpha >0$. \\ \\
(ii) The function $f$ is in $C^k(\R_{\geq 0})$. For $\alpha \in (0,1/2)$ we have $|f^{(j)}(x)|\leq c_j |x|^{-j}$ for $x\geq 1$,
while for $\alpha \geq 1/2$ there exists a $v>0$ such that $\alpha -1/2<v$ and $|f^{(j)}(x)|\leq c_j |x|^{-j-v}$ for $x\geq 1$ ($0\leq j\leq k$). \\ \\
(iii) $L$ is a symmetric L\'evy processes with L\'evy measure $\nu$ and without Gaussian part. L\'evy measure $\nu$ satisfies
\begin{align*}
&\int_{|x|\geq 1} |x|^{1/(1-\alpha)} \nu(dx)<\infty  \qquad \text{for } \alpha \in (0,1/2), \\
&\int_{|x|\geq 1} |x|^{1/(v+1-\alpha)} \nu(dx)<\infty  \qquad \text{for } \alpha\geq 1/2.  
\end{align*}
Another important parameter in our limit theory is the \textit{Blumenthal-Getoor index} of the driving L\'evy motion 
$L$, which is defined by 

\begin{align} \label{bgindex}
\beta := \inf \left\{r\geq 0:~ \sum_{s\in [0,1]} |\Delta L_s|^r<\infty  \right\}, \qquad \Delta L_s = L_s - L_{s-}. 
\end{align}
Obviously, finite activity L\'evy processes have Blumenthal-Getoor index $\beta =0$ while L\'evy processes with finite
variation satisfy $\beta \leq 1$. In general, the Blumenthal-Getoor index is a non-random number $\beta \in [0,2]$
and it can characterized by the L\'evy measure $\nu$ of $L$ as follows:
\[
\beta = \inf \left\{r\geq 0:~ \int_{-1}^1 |x|^r \nu (dx)<\infty  \right\}.
\]
The latter implies that $\beta$-stable L\'evy processes with $\beta \in (0,2)$ have Blumenthal-Getoor index $\beta$.
  
We recall the definition of $k$th order increments $\Delta_{i,k}^n Y$ introduced at \eqref{filter}
and consider the power variation 
\[
V(Y,p,k;\Delta_n)_t = \sum_{i=k}^{[t/\Delta_n]} |\Delta_{i,k}^{n} Y|^p.
\]
The following theorem from \cite{BLP} determines the first order structure of the statistic $V(Y,p,k;\Delta_n)_t$.
Notice that the results below are stated for a fixed $t>0$.

\begin{theo} \label{th7}
Assume that condition (A') holds and fix $t>0$.
\begin{itemize}
\item[(i)] If $\alpha \in (0,k-1/p)$ and $p>\beta$, we obtain the stable convergence
\begin{align} \label{part1}
\Delta_n^{-\alpha p} V(Y,p,k;\Delta_n)_t \stab |f(0)|^p \sum_{m:T_m\in [0,t]} |\Delta L_{T_m}|^p \left( \sum_{l=0}^{\infty} |
h(l+U_m)|^p \right),    
\end{align} 
where $(U_m)_{m\geq 1}$ are i.i.d. $\mathcal{U}([0,1])$-distributed random variables independent of the original
$\sigma$-algebra $\mathcal{F}$ and the function $h=h_0$ is defined at \eqref{hdef}.
\item[(ii)] Assume that $L$ is a symmetric $\beta$-stable L\'evy process with $\beta \in (0,2)$. 
If $\alpha \in (0,k-1/\beta)$ and $p<\beta$
then it holds
\begin{align} \label{part2}
\Delta_n^{1-p(\alpha + 1/\beta)} V(Y,p,k;\Delta_n)_t \toop t c_p, \qquad c_p:=\E[|L_1(k)|^p]<\infty ,
\end{align}
where the process $L_t(k)$ is defined as
\[
L_t(k):= \int_{\R} h(s) dL_s
\]
and the function $h=h_0$ is defined at \eqref{hdef}.
\item[(iii)] When $\alpha > k-1/p,~p>\beta$ or $\beta > k-1/\beta,~ p<\beta$, we deduce
\begin{align} \label{part3}
\Delta_n^{1-pk}V(Y,p,k;\Delta_n)_t \toop \int_0^t |F_k(u)|^p du , \qquad  F_k(u):= \int_{-\infty}^u g^{(k)}(u-s) dL_s. 
\end{align}
%\item[(iv)] (Critical case I) If $\alpha =k-1/p$ and $p>\beta$, we conclude that 
%\begin{align} \label{part4}
%\frac{\Delta_n^{-\alpha p}}{\log n}V(Y,p,k;\Delta_n)_t \toop |f(0)|^p \sum_{m:T_m\in [0,t]} |\Delta L_{T_m}|^p. 
%\end{align}
%\item[(v)] (Critical case II) If $\alpha = k-1/\beta$ and $p<\beta$, we obtain 
%\begin{align} \label{part5}
%\frac{\Delta_n^{1-p(\alpha + 1/\beta)}}{(-\log \Delta_n)^{p/\beta }} V(Y,p,k;\Delta_n)_t  \toop tc_p.
%\end{align}
\end{itemize}
\end{theo}
We remark that the result of Theorem \ref{th7}(i) is sharp in a sense that the conditions $\alpha \in (0,k-1/p)$ and $p>\beta$
are sufficient and (essentially) necessary to conclude \eqref{part1}. Indeed, since $|h(l+U_m)|\leq \text{const.}l^{\alpha -k}$ for $l\geq 1$, 
we obtain that 
\[
\sum_{l=0}^{\infty} |h(l+U_m)|^p \leq \text{const}<\infty  
\]
when $\alpha \in (0,k-1/p)$, and on the other hand $\sum_{m:T_m\in [0,t]} |\Delta L_{T_m}|^p<\infty$ for $p>\beta$, which follows from the definition
of the Blumenthal-Getoor index.

The idea behind the proof of Theorem \ref{th7}(ii) has been described in Section \ref{sec2.3}. Note that for $k=1$ the random variable
$L_1(1)$ coincides with the increments  $\widetilde X_1- \widetilde X_0$ of the fractional $\beta$-stable L\'evy motion introduced in \eqref{flm1}.
Following the mathematical intuition of the aforementioned discussion, the limit in \eqref{part2} is not really surprising. Also notice
that the limit is indeed finite since $p<\beta$.

We remark that for values of $\alpha$ close to $k-1/p$ or $k-1/\beta $ in Theorem \ref{th7}(iii), the function $g^{(k)}$ explodes at $0$. This leads
to unboundedness of the process $F_k$ defined at \eqref{part3}. Nevertheless, under conditions of Theorem \ref{th7}(iii), the limiting process
is still finite.

\begin{rem} \label{rem14} \rm
Notice that the critical cases, i.e. $\alpha =k-1/p, p>\beta $ and $\alpha =k-1/, p<\beta $, are not described in Theorem \ref{th7}. In this cases
an additional log factor appears, and for $\alpha =k-1/p, p>\beta $ the mode of convergence changes from stable convergence to convergence in probability
(clearly the limits change too). We refer to \cite{BLP} for a detailed discussion of critical cases. \qed
\end{rem}    

\begin{rem} \label{rem14a} \rm
The asymptotic results of Theorem \ref{th7} uniquely identify the parameters $\alpha$ and $\beta$. First of all, note that the convergence
rates in \eqref{part1}-\eqref{part3} are all different under the corresponding conditions. Indeed, it holds that
\[
p(\alpha +1/\beta)-1<\alpha p< pk -1,
\]
since in case (i) we have $\alpha <k-1/p$ and in case (ii) we have $p<\beta$. Hence, computing the statistic $V(Y,p,k;\Delta_n)_t$ at log scale
for all $p\in [0,2]$  identifies the parameters $\alpha$ and $\beta$. \qed   
\end{rem}

\begin{rem} \label{rem15} \rm
A related study of the asymptotic theory is presented in \cite{BCI},  who investigated the fine structure of L\'evy moving average processes driven by
a truncated $\beta$-stable L\'evy motion. The authors showed the result of Theorem \ref{th7}(ii) (see Theorem 5.1 therein), whose prove was 
however incorrect, since it was based on the computation of the variance that diverges to infinity. 

In a recent work \cite{Gl} extended the law of large numbers in Theorem \ref{th7}(ii) to integral processes driven by fractional L\'evy
motion. The main idea relies on the mathematical intuition described in Remark \ref{rem8}. \qed      
\end{rem}  
The next theorem demonstrates a central limit theorem associated with Theorem \ref{th7}(ii) (see \cite{BLP}).

\begin{theo} \label{th8}
Assume that condition (A') holds and fix $t>0$. Let $L$ be a symmetric $\beta$-stable L\'evy process with characteristic triplet
$(0,0,c|x|^{-1-\beta}dx)$ and $\beta \in (0,2)$. 
When $k\geq 2$, $\alpha \in (0,k-2/\beta)$ and $p<\beta/2$
then it holds
\begin{align} \label{newclt}
\Delta_n^{-1/2} \left( \Delta_n^{1-p(\alpha + 1/\beta)} V(Y,p,k;\Delta_n)_t - t c_p \right) \schw \mathcal N(0, t\eta^2), 
\end{align}
where the quantity $\eta^2$ is defined via
\begin{align*}
\eta^2 &= \theta(0) + 2 \sum_{i=1}^\infty \theta(i)  , \\[1.5 ex]
\theta(i) &= a_p^{-2} \int_{\R^2} \frac{1}{|s_1s_2|^{1+p}} \psi_i(s_1,s_2) ds_1 ds_2, \\[1.5 ex]
\psi_i(s_1,s_2) &= \exp \left(-c|f(0)|^{\beta} \int_{\R} |s_1 h(x) -s_2h(x+i)|^{\beta } dx \right),  \\[1.5 ex]
&- \exp \left(-c|f(0)|^{\beta} \int_{\R} |s_1 h(x)|^{\beta} + |s_2h(x+i)|^{\beta} dx \right),
\end{align*}
where the function $h$ is defined at \eqref{hdef} and $a_p:= \int_{\R} (\exp(iu) -1) |u|^{-1-p} du$.
\end{theo}
Let us explain the various conditions of Theorem \ref{th8}. The assumption $p<\beta/2$ ensures the existence of variance of the statistic $V(Y,p,k;\Delta_n)_t$.
The validity range of the central limit theorem ($\alpha \in (0,k-2/\beta)$) is smaller than the validity range of the law of large numbers
in Theorem \ref{th7}(ii) ($\alpha \in (0,k-1/\beta)$). It is not clear which limit distribution appears in case $\alpha \in (k-2/\beta, k-1/\beta)$.
A more severe assumption is $k\geq 2$, which excludes the first order increments. The limit theory in this case is also unknown.

\begin{rem} \label{rem17} \rm
Let us explain a somewhat complex form of the variance $\eta^2$. A major problem of proving Theorems \ref{th7}(ii) and \ref{th8}
is that neither the expectation of $|\Delta_{i,k}^{n} Y|^p$ nor its variance can be computed directly. However, the identity
\[
|x|^p = a_p^{-1} \int_{\R} (\exp(iux) -1) |u|^{-1-p} du \qquad \text{for } p\in (0,1),
\] 
which can be shown by substitution $y=ux$, turns out to be a useful instrument. Indeed, for any deterministic function $\varphi :\R\rightarrow \R$
satisfying the conditions of Remark \ref{rem1}, it holds that
\[
\E\left[ \exp \left(iu \int_{\R} \varphi_s dL_s \right)\right] = \exp \left(-c|u|^{\beta} \int_{\R} |\varphi_s|^{\beta} ds \right).
\] 
This two identities are used to compute the variance of the statistic $V(Y,p,k;\Delta_n)_t$ and they are both reflected in the formula for
the quantity $\theta(i)$.
\end{rem}

\begin{rem} \label{rem16} \rm
As in the case of a Gaussian driver, Theorem \ref{th7}(ii) might be useful for statistical applications. Indeed, for any fixed $t>0$, it holds
that
\[
S_n=\frac{V(Y,p,k;2\Delta_n)_t}{V(Y,p,k;\Delta_n)_t} \toop 2^{p(\alpha +1/\beta)}, 
\] 
under conditions of Theorem \ref{th7}(ii). Thus, a consistent estimator of $\alpha$ (resp. $\beta$) can be constructed given the knowledge
of $\beta$ (resp. $\alpha$) and the validity of conditions $\alpha \in (0,k-1/\beta)$ and $p<\beta$. A bivariate version of the central limit theorem
\eqref{newclt} for frequencies $\Delta_n$ and $2\Delta_n$ would give a possibility to construct feasible confidence regions. \qed     
\end{rem}    
Obviously, the presented limiting results still need to be extended to spatio-temporal setting. Asymptotic theory for ambit fields
observed on a grid in time-space or along a curve would be extremely useful for statistical analysis of turbulent flows.


\begin{thebibliography}{1}

\bibitem{AE78} D.J. Aldous and G.K. Eagleson (1978):  On mixing and stability of limit theorems.
\textit{Annals of Probability} 6(2), 325--331.

\bibitem{AMN} E. Alos, O. Mazet and D. Nualart (2001): Stochastic calculus with respect to Gaussian
processes. \emph{Annals of Probability} 29(2), 766--801.

\bibitem{BBPV} O.E. Barndorff-Nielsen, F.E. Benth, J. Pedersen and A. Veraart (2014): On stochastic integration for volatility modulated L\'evy–driven
Volterra processes. \emph{Stochastic Processes and Their Applications} 124(1), 812-847.

\bibitem{BBS} O.E. Barndorff-Nielsen, F.E. Benth and B. Szozda (2013): On stochastic integration for volatility modulated
Brownian-driven Volterra processes via white noise analysis. \emph{Working paper}, available at arXiv:1303.4625.

\bibitem{BBC} 
O.E. Barndorff-Nielsen, F.E. Benth and A. Veraart (2011): Modelling electricity forward markets by ambit fields. \textit{Working paper}. Available at SSRN:http://ssrn.com/abstract=1938704.

\bibitem{BBC2} 
O.E. Barndorff-Nielsen, F.E. Benth and A. Veraart (2012): Recent advances in ambit stochastics with a view towards
tempo-spatial stochastic volatility/intermittency. \emph{Working paper.} Available at arXiv:1210.1354.

\bibitem{BCP09} O.E. Barndorff-Nielsen,  J.M. Corcuera and M. Podolskij (2009):
Power variation for Gaussian processes with stationary increments. 
{\em Stochastic Processes and Their Applications} 119, 1845--865.
 
\bibitem{BCP11} O.E. Barndorff-Nielsen,  J.M. Corcuera and M. Podolskij (2011): 
Multipower variation for Brownian semistationary processes. 
\textit{Bernoulli} 17(4), 1159--1194. 

\bibitem{BCP13} O.E. Barndorff-Nielsen,  J.M. Corcuera and M. Podolskij (2013): Limit theorems for functionals of higher order 
differences of Brownian semi-stationary processes. In {\em Prokhorov and Contemporary Probability Theory: 
In Honor of Yuri V. Prokhorov}, eds. A.N. Shiryaev, S.R.S. Varadhan and E.L. Presman. Springer.

\bibitem{BJJS07}  O.E. Barndorff-Nielsen, E.B.V. Jensen, K.Y. J\'onsd\'ottir and J. Schmiegel (2007): 
Spatio-temporal modelling - with a view to biological growth. In B. Finkenstädt, L. Held and V. Isham: 
\emph{Statistical Methods for Spatio-Temporal Systems}. London: Chapman and Hall/CRC, 47-75.

\bibitem{BGJPS06} O.E. Barndorff-Nielsen,  S.E. Graversen, 
J. Jacod, M. Podolskij and N. Shephard (2006): A central limit theorem
for realised power and bipower variations of continuous
semimartingales. In: Kabanov, Yu., Liptser, R., Stoyanov, J.~(eds.),
From Stochastic Calculus to Mathematical Finance. Festschrift in
Honour of A.N.~Shiryaev, pp.~33--68. Springer, Heidelberg.

\bibitem{BPS} O.E. Barndorff-Nielsen, M. Pakkanen and J. Schmiegel (2013): 
Assessing relative volatility/intermittency/energy dissipation. \emph{Working paper}, available at arXiv:1304.6683.

\bibitem{BS07} O.E. Barndorff-Nielsen and J. Schmiegel (2007): Ambit
processes; with applications to turbulence and cancer growth. In F.E. Benth,
Nunno, G.D., Linstr\o m, T., \O ksendal, B. and Zhang, T. (Eds.):\textrm{\
\emph{Stochastic Analysis and Applications: The Abel Symposium 2005}}.
Heidelberg: Springer. Pp. 93--124.\textrm{\ }

\bibitem{BS08a} O.E. Barndorff-Nielsen and J. Schmiegel (2008): Time
change, volatility and turbulence. In A. Sarychev, A. Shiryaev, M. Guerra
and M.d.R. Grossinho (Eds.): \emph{Proceedings of the Workshop on
Mathematical Control Theory and Finance}. Lisbon 2007. Berlin: Springer. Pp.
29--53.

\bibitem{BS09} O.E. Barndorff-Nielsen and J. Schmiegel (2009): Brownian
semistationary processes and volatility/intermittency. In: Albrecher, H.,
Runggaldier, W., Schachermayer, W. (Eds.): \emph{Advanced Financial Modelling%
}, 1--26, Germany: Walter de Gruyter.

\bibitem{B} A. Basse (2008): Gaussian moving averages and semimartingales. \emph{Electronic Journal of Probability} 
13(39), 1140--1165.

\bibitem{BLP} A. Basse-O'Connor, R. Lachieze-Rey and M. Podolskij (2014): Limit theorems for L\'evy moving average
processes. \emph{Working paper}. 

\bibitem{BP} A. Basse and J. Pedersen (2009): L\'evy driven moving averages and semimartingales. 
\emph{Stochastic Processes and Their Applications} 119(9), 2970--2991.

\bibitem{BR} A. Basse-O'Connor and J. Rosinski (2012): Structure of infinitely divisible semimartingales.
\emph{Working paper}, available at arXiv:1209.1644v2.

\bibitem{BCI} A. Benassi, S. Cohen and J. Istas (2004): On roughness indices for fractional fields.
\emph{Bernoulli} 10(2), 357--373.

\bibitem{BLS} C. Bender, A. Lindner and M. Schicks (2012): Finite variation of fractional Le´vy processes. 
\emph{Journal of Theoretical Probability} 25(2), 594--612.

\bibitem{BEV} F.E. Benth, H. Eyjolfsson and  A.E.D. Veraart (2013): Approximating L\'evy semistationary processes via
Fourier methods in the context of power markets. \emph{Working paper}.

\bibitem{BS13} F.E. Benth and A. S\"u{\ss} (2013): Integration theory for infinite dimensional volatility modulated Volterra processes.
\emph{Working paper}, available at arXiv:1303.7143. 

\bibitem{BJ} K. Bichteler and J. Jacod (1983): Random measures and stochastic integration. 
\emph{Theory and Application of Random Fields. 
Lecture Notes in Control and Information Sciences} 49,  1--18. 

\bibitem{BM83} P. Breuer and P. Major (1983): Central limit theorems for nonlinear functionals of
Gaussian fields. \textit{Journal of Multivariate Analysis} 13(3), 425--441.

\bibitem{CK} C. Chong and C. Kl\"uppelberg (2013): Integrability conditions for space-time stochastic integrals: theory
and applications. \emph{Working paper}, available at arXiv:1303.2468. 

\bibitem{CHPP13} J.M. Corcuera, E. Hedevang, M. Pakkanen and M. Podolskij (2013): Asymptotic theory for Brownian semi-stationary processes with application to turbulence. \emph{Stochastic Processes and Their Applications} 123, 2552-2574. 

\bibitem{CNP14}  J.M. Corcuera, D. Nualart and M. Podolskij (2014): Asymptotics of weighted random sums. \emph{Working paper},
available at http://arxiv.org/abs/1402.1414. 

\bibitem{CNW} J. M. Corcuera,  D. Nualart and J.H.C. Woerner (2006): 
Power variation of some integral fractional processes. \textit{Bernoulli} 12(4), 713--735.

\bibitem{DQ} R.C. Dalang and L. Quer-Sardanyons (2011): Stochastic integrals for spde's: a comparison.
\emph{Expositiones Mathematicae} 12(1), 67-109.

\bibitem{NOP} G. Di Nunno, B. Oksendal and F. Proske (2009): \emph{Malliavin calculus for L\'evy processes with
applications to finance}. Springer-Verlag, Berlin.

\bibitem{GP} K. G\"artner and M. Podolskij (2014): On non-standard limits of Brownian semi-stationary processes.
\textit{Working paper.} Available at arXiv:1403.6484.

\bibitem{Gl} S. Glaser (2014): A law of large numbers for the power variation of fractional L\'evy processes. \emph{Working paper}. 

\bibitem{J08}
J. Jacod (2008): Asymptotic properties of realized power variations and
related functionals of semimartingales. {\em Stoch. Process. Appl.}
\textbf{118}, 517--559.

\bibitem{JP} J. Jacod and P.E. Protter (2012): \emph{Discretization of processes.} Springer.  

\bibitem{JS02} J. Jacod and A.N. Shiryaev (2002):
{\em Limit theorems for stochastic
  processes}, 2d Edition. Springer Verlag: Berlin.
  
\bibitem{JY} T. Jeulin and M. Yor (1993): Moyennes mobiles et semimartingales. S\'eminaire de Probabilit\'es 
XXVII (1557), 53--77. 
  
  
\bibitem{K}  F.B. Knight (1992): \emph{Foundations of the Prediction Process}. Volume 1 of Oxford Studies in Probability. New York: The Clarendon Press Oxford University Press. Oxford Science Publications.

\bibitem{M} T. Marquardt (2006): Fractional L\'evy processes with an application to long memory moving
average processes. \emph{Bernoulli} 12, 1090–-1126.

\bibitem{N} D. Nualart (2006): \textit{The Malliavin calculus and
related topics.} (2nd edition). Springer, Berlin.

\bibitem{NP05} D. Nualart and G. Peccati (2005): Central limit theorems for multiple stochastic integrals. \textit{Ann. Probab.} 33(1), 177--193.

\bibitem{Pa} M. Pakkanen (2013): Limit theorems for power variations of ambit fields driven by white noise. \textit{Working paper},
available at arXiv:1301.2107.

\bibitem{PR} M. Pakkanen and A. R\'eveillac (2014): Functional limit theorems for generalized variations of the fractional
Brownian sheet. \emph{Working paper.} Available at arXiv:1404.2822. 

\bibitem{PT05} G. Peccati and C.A. Tudor (2005): Gaussian limits for vector-values multiple stochastic integrals. In: \emph{S\'eminaire de Probabilit\'es XXXVIII}, Lecture Notes in Math., 1857. Berlin: Springer. Pp. 247--193.  

\bibitem{Pr} P. Protter (1985): Volterra equations driven by semimartingales. \emph{Annals of Probability} 13(2), 519--530. 

\bibitem{RR} B. Rajput and J. Rosinski(1989): Spectral representation of infinitely divisible distributions. {\em Probability
Theory and Related Fields} 82, 451--487.
  
\bibitem{R63} A. R\'enyi (1963): On stable sequences of events. \textit{Sankhy\=a Ser. A} 25, 293--302.

\bibitem{Re} A. R\'eveillac (2009): Estimation of quadratic variation for two-parameter diffusions.
\textit{Stochastic Processes and Their Applications} 119(5), 1652--1672.

\bibitem{T} M. Taqqu (1979): Convergence of integrated processes of arbitrary Hermite rank. 
\textit{Z. Wahrsch. Verw. Gebiete} 50(1), 53--83. 

\bibitem{W} J. Walsch (1986): An introduction to stochastic partial differential equations. 
{\em \'Ecole d'Et\'e de Prob. de St. Flour XIV, Lect. Notes in Math.} 1180,  265--439.




\end{thebibliography}
\end{document}